\documentclass[12pt,reqno]{amsart}
\usepackage[notref,notcite]{}
\usepackage{graphicx}
\usepackage{amssymb}
\usepackage{amsmath}

\numberwithin{equation}{section}

\newtheorem{theorem}{Theorem}[section]
\newtheorem{lemma}[theorem]{Lemma}
\newtheorem{proposition}[theorem]{Proposition}

\def\l{\langle}

\def\C{\mathbb C}

\def\R{\mathbb R}
\def\N{\mathbb N}
\def\Z{\mathbb Z}

\newcommand{\cout}[1]{}

\begin{document}
\title[Rank 3 Killing tensor fields on a 2-torus]{Killing tensor fields of third rank on a two-dimensional Riemannian torus}

\author{Vladimir A. Sharafutdinov}
\thanks{Partially supported by RFBR Grant 20-51-15004 (joint French--Russian grant).}
\address{Sobolev Institute of Mathematics\\
4 Koptyug avenue, Novosibirsk, 630090, Russia}
\email{sharaf@math.nsc.ru}


%\subjclass[2000]{Primary 35R30; Secondary 35P99}

\begin{abstract}
A rank $m$ symmetric tensor field on a Riemannian manifold is called a Killing field if the symmetric part of its covariant derivative is equal to zero.
Such a field determines the first integral of the geodesic flow which is a degree $m$ homogeneous polynomial in velocities. There exist global isothermal coordinates on a two-dimensional Riemannian torus such that the metric is of the form $ds^2=\lambda(z)|dz|^2$ in the coordinates. The torus admits a third rank Killing tensor field if and only if the function $\lambda$ satisfies the equation $\Re\big(\frac{\partial}{\partial z}\big(\lambda(c\Delta^{-1}\lambda_{zz}+a)\big)\big)=0$ with some complex constants $a$ and $c\neq0$. The latter equation is equivalent to some system of quadratic equations relating Fourier coefficients of the function $\lambda$. If the functions $\lambda$ and $\lambda+\lambda_0$ satisfy the equation for a real constant $\lambda_0\neq0$, then there exists a non-zero Killing vector field on the torus.
\end{abstract}

\maketitle

\section{Introduction}

We first recall the definition of a Killing tensor field on an arbitrary Riemannian manifold.

Given a Riemannian manifold $(M,g)$, let $S^m\tau'_M$ be the bundle of rank $m$ symmetric tensors. The latter notation will mostly be abbreviated to $S^m$ assuming the manifold under consideration to be known from the context. The space $C^\infty(S^m)$ of smooth sections of the bundle consists of degree $m$ symmetric tensor fields. The sum
$
S^*=\bigoplus_{m=0}^\infty S^m
$
is a bundle of graded commutative algebras with respect to the product $fh=\sigma(f\otimes h)$, where $\sigma$ is the symmetrization.

The differential operator
$
d=\sigma\nabla:C^\infty(S^m)\rightarrow C^\infty(S^{m+1}),
$
where $\nabla$ is the covariant derivative with respect to the Levi-Chivita connection, is called the {\it inner derivative}. A tensor field $f\in C^\infty(S^m)$ is said to be a {\it Killing tensor field} if
\begin{equation}
df=0.
                  \label{1.1}
\end{equation}
In the case of $m=1$ we say on a {\it Killing vector field}.
The operator $d$ is related to the product by the Leibnitz formula $d(fh)=(df)h+f(dh)$ that implies the statement: if $f$ and $h$ are Killing tensor fields, then $fh$ is also a Killing field. A Killing tensor field $f\in C^\infty(S^m)\ (m\neq2)$ is said to be {\it irreducible}, if it cannot be represented as a finite sum $f=\sum_i u_iv_i$, where all $u_i$ and $v_i$ are Killing fields of positive ranks. In the case of $m=2$ we additionally require $f$ to be different of
$cg\ (c=\mbox{const})$. The requirement eliminates the metric tensor from irreducible Killing fields.

Being written in coordinates for a rank $m$ tensor field, \eqref{1.1} is a system of ${n+m}\choose {m+1}$ linear first order differential equations in ${n+m-1}\choose {m}$ coordinates of $f$, where $n=\mbox{dim}M$. Since the system is overdetermined, not every Riemannian manifold admits nonzero Killing tensor fields. The two-dimensional case is of the most interest since the degree of the overdetermination is equal to 1 in this case. In the two-dimensional case we obtain one equation on the metric $g$ after eliminating all coordinates of $f$ from the system \eqref{1.1}. The possibility of such elimination is rather problematic in the general case, but it will be realized below in the case of $m=3$ for a 2-torus.

Let $\pi:TM\rightarrow M$ be the tangent bundle. We denote points of the manifold $TM$ by pairs $(x,\xi)$, where $x\in M$ and $\xi\in T_xM$. If $(U;x^1,\dots,x^n)$ is a local coordinate system on $M$ with the domain $U\subset M$, then the corresponding local coordinate system $(\pi^{-1}(U);x^1,\dots,x^n,\xi^1,\dots,\xi^n)$ is defined on $TM$, where $\xi=\xi^i\frac{\partial}{\partial x^i}$. A tensor field $f\in C^\infty(S^m)$ is written as $f=f_{i_1\dots i_m}\,dx^{i_1}\dots dx^{i_m}$ in local coordinates. Given a tensor field $f\in C^\infty(S^m)$, we define the function $F\in C^\infty(TM)$ by $F(x,\xi)=f_{i_1\dots i_m}(x)\,\xi^{i_1}\dots\xi^{i_m}$. The map $f\mapsto F$ is independent of the choice of local coordinates and identifies the algebra $C^\infty(S^*)$ with the subalgebra of $C^\infty(TM)$ which consists of functions polynomially depending on $\xi$.

Let $H$ be the vector field on $TM$ generating the geodesic flow.
Let $\Omega M\subset TM$ be the manifold of unit tangent vectors. Since the geodesic flow preserves the norm of a vector, $H$ can be considered as a first order differential operator on $\Omega M$. The operators $d$ and $H$ are related as follows: if $f\in C^\infty(S^m)$ and $F=f_{i_1\dots i_m}\,\xi^{i_1}\dots\xi^{i_m}\in C^\infty(TM)$ is the corresponding polynomial, then
$HF=(df)_{i_1\dots i_{m+1}}\,\xi^{i_{1}}\dots\xi^{i_{m+1}}$. In particular, $f$ is a Killing tensor field if and only if $HF=0$, i.e., if $F$ is a first integral for the geodesic flow. Thus, the problem of finding Killing tensor fields is equivalent to the problem of finding first integrals of the geodesic flow which polynomially depend on $\xi$.

\bigskip

We proceed to considering Killing tensor fields on the 2-torus. The following question remains open although it has been considered in many works \cite{KD,Sh,BM,MS}.

($\ast$) {\it Does there exist a Riemannian metric on the 2-torus which admits an irreducible Killing tensor field of rank $m\ge3$?}

The present paper is devoted to the question for $m=3$. The question remains open, our results are of a very particular character. Nevertheless we hope the paper will serve for a further progress in this hard problem.

Recall \cite[\S6.5]{BF} that there exists a global isothermal coordinate system on the two-dimensional torus ${\mathbb T}^2$ furnished with a Riemannian metric $g$. More precisely, there exists a lattice $\Gamma\subset{\mathbb R}^2={\mathbb C}$ such that ${\mathbb T}^2={\mathbb C}/\Gamma$ and the metric $g$ is expressed by the formula
\begin{equation}
g=\lambda(x,y)(dx^2+dy^2)=\lambda(z)|dz|^2\quad(z=x+iy),
                            \label{1.2}
\end{equation}
where $\lambda(z)$ is a $\Gamma$-periodic smooth positive function on the plane. Global isothermal coordinates are defined up to changes of the kind $z=az'+b$ or $z=a\bar z{}'+b$ with complex constants $a\neq 0$ and $b$. We will widely use changes $z=az'\ (0\neq a\in\C)$ that geometrically mean the possibility to rotate the lattice $\Gamma$ through an arbitrary angle around the origin and to stretch (squeeze) the lattice with respect to the origin with an arbitrary positive coefficient.

The question ($\ast$) is completely investigated in the cases of $m=1$ and of $m=2$. A Riemannian 2-torus admits a non-trivial (not identically equal to zero) Killing vector field if and only if the metric is of the form \eqref{1.2} in some global isothermal coordinate system, where
\begin{equation}
\lambda(x,y)=\mu(x).
                            \label{1.3}
\end{equation}
The torus admits a second rank irreducible Killing tensor field if and only if the metric is of the form \eqref{1.2} in some global isothermal coordinate system, where
\begin{equation}
\lambda(x,y)=\mu(x)+\nu(y)
                            \label{1.4}
\end{equation}
and both $\mu$ and $\nu$ are non-constant functions.
Such metrics are named {\it Liouville metrics}.

We return to considering an arbitrary Riemannian torus $({\C}/\Gamma,\lambda|dz|^2)$. Since $\lambda$ is a smooth $\Gamma$-periodic function, it can be represented by the Fourier series
\begin{equation}
\lambda(x,y)=\sum\limits_{(n_1,n_2)\in\Gamma'}{\hat\lambda}_n\,e^{i(n_1x+n_2y)}
                         \label{1.5}
\end{equation}
with coefficients rapidly decaying as $|n_1|+|n_2|\rightarrow\infty$. Here
$$
\Gamma'=\{k=(k_1,k_2)\in{\R}^2\mid k_1n_1+k_2n_2\in 2\pi{\Z}\ \mbox{for all}\ n=(n_1,n_2)\in\Gamma\}
$$
is the {\it dual lattice} for $\Gamma$. The following statement is one of results of the present work.

\begin{theorem} \label{Th1.1}
A Riemannian torus $({\C}/\Gamma,\lambda|dz|^2)$ admits a non-trivial rank 3 Killing tensor field if and only if Fourier coefficients of the function $\lambda$ satisfy the equations
\begin{equation}
\sum\limits_{k\in\Gamma'\setminus\{0\}}
\!\!\!\!\frac{c_1(\!-\!n_1k_1^2\!+\!2n_2k_1k_2\!+\!n_1k_2^2)+c_2(\!-\!n_2k_1^2\!-\!2n_1k_1k_2\!+\!n_2k_2^2)}{k_1^2+k_2^2}
{\hat\lambda}_k{\hat\lambda}_{n-k}
=(a_1n_1+a_2n_2){\hat\lambda}_n
                         \label{1.6}
\end{equation}
for all $n\in\Gamma'$ with some real constants $(a_1,a_2)$ and $(c_1,c_2)\neq(0,0)$.
\end{theorem}

The system \eqref{1.6} can be written as one pseudodifferential equation. To this end we observe that the inverse operator for the Laplacian
$\Delta=\frac{\partial^2}{\partial x^2}+\frac{\partial^2}{\partial y^2}$ can be defined on the space of $\Gamma$-periodic functions by
\begin{equation}
\Delta^{-1}e^{i(n_1x+n_2y)}=-\frac{1}{n_1^2+n_2^2}\,e^{i(n_1x+n_2y)}\ (0\neq n\in\Gamma'),\quad
\Delta^{-1}\,\mbox{const}=0.
                         \label{1.7}
\end{equation}
For a real $\Gamma$-periodic function $\lambda\in C^2(\R^2)$, the system \eqref{1.6} is equivalent to the equation
\begin{equation}
\frac{\partial}{\partial z}\big(\lambda(c\Delta^{-1}\lambda_{zz}+a)\big)
+\frac{\partial}{\partial\bar z}\big(\lambda(\bar c\Delta^{-1}\lambda_{\bar z\bar z}+\bar a)\big)=0,
                         \label{1.8}
\end{equation}
where $a=\frac{1}{4}(a_1+ia_2)$ and $c=c_1+ic_2\neq0$. The equation \eqref{1.8} can be considered as the result of the above-mentioned elimination of coordinates of the field $f$ from the system \eqref{1.1}.

For every real function $\mu\in C^\infty(\R)$ and every pair of real constants $(\alpha,\beta)$, the function
\begin{equation}
\lambda=\mu(\alpha x+\beta y)
                         \label{1.9}
\end{equation}
solves the equation \eqref{1.8} with respectively chosen constants $a$ and $c$. One can easily find conditions on the function $\mu$ and coefficients $(\alpha,\beta)$ which guarantee $\Gamma$-periodicity of the function \eqref{1.9}. Such solutions will be called {\it one-dimensional solutions}. The geometric sense of such solutions is obvious: If the equation \eqref{1.8} possesses a $\Gamma$-periodic one-dimensional solution, then the metric can be transformed to the form \eqref{1.2}--\eqref{1.3} by some rotation of the lattice $\Gamma$ and the torus admits a non-trivial Killing vector field $f$. In such a case
$f^3=f\otimes f\otimes f$ is the reducible Killing tensor field of third rank. We are thus interested in looking for solutions to the equation \eqref{1.8} (or system {1.6}) which are not one-dimensional ones.

Following M.V. Denisova and V.V. Kozlov \cite{KD}, we define the {\it spectrum} of a $\Gamma$-periodic function $\lambda\in C^\infty(\R^2)$ as the set of nodes $n$ of the lattice $\Gamma'$ such that ${\hat\lambda}_n\neq0$. A $\Gamma$-periodic solution to the equation \eqref{1.8} is a one-dimensional solution if and only if its spectrum is contained in the intersection $L\cap\Gamma'$, where $L\subset\R^2$ is a line through the origin.

For $m=3$ the question ($\ast$) is equivalent to the following one:

($\ast\ast$) {\it Does there exist a lattice $\Gamma\subset\R^2$ such that the equation \eqref{1.8} (or system \eqref{1.6}) has a $\Gamma$-periodic positive solution $\lambda\in C^\infty(\R^2)$ which is not a one-dimensional solution?}

The system \eqref{1.6} involves four real constants $(a_1,a_2,c_1,c_2)$. It is clear that only three pairwise ratios $(a_1:a_2:c_1:c_2)$ are essential. We emphasize that a lattice $\Gamma'$ is also considered as an unknown in \eqref{1.6}. Up to rotation and homothety, a lattice is determined by two real constants. Thus the system \eqref{1.6} involves five real constants that should be considered as unknowns as well as Fourier coefficients ${\hat\lambda}_k$.

Reality of the function $\lambda$ is guaranteed by the following {\it parity requirement} for a solution to the system \eqref{1.6}:
\begin{equation}
{\hat\lambda}_{-k}=\overline{{\hat\lambda}_k}\quad (k\in\Gamma').
                         \label{1.10}
\end{equation}
Smoothness of $\lambda$ is guaranteed by the {\it decay condition}: for every $N\in\N$
\begin{equation}
|{\hat\lambda}_{k}|\le C_N(k_1^2+k_2^2+1)^{-N}.
                         \label{1.11}
\end{equation}

Most probably, positiveness of $\lambda$ is the most difficult requirement in the question ($\ast\ast$). The Bochner -- Hinchin theorem, that is popular in probability theory \cite{B}, allows to state the requirement in terms of Fourier coefficients ${\hat\lambda}_n$. Unfortunately, the statement of
Bochner -- Hinchin's theorem is not connected to the system \eqref{1.6} (at least the author does not see such a connection). Therefore the following weaker version of the question ($\ast\ast$) seems to be reasonable:

($\ast\ast\ast$) {\it Does there exist a lattice $\Gamma\subset\R^2$ such that the equation \eqref{1.8} (or system \eqref{1.6}) has a $\Gamma$-periodic real solution $\lambda\in C^\infty(\R^2)$ which is not a one-dimensional solution and satisfies ${\hat\lambda}_0>0$?}

In the cases of $m=1$ and of $m=2$, our problem has the following specifics. If a function $\lambda$ is either of the form \eqref{1.3} or of the form \eqref{1.4}, then $\lambda+\lambda_0$ is of the same form for any real constant $\lambda_0$. Do the specifics preserve for $m=3$? If the answer to the question was positive, then the main difficulty of our problem (the positiveness requirement) would disappear. Indeed, if $\lambda\in C^\infty(\R^2)$ is a real $\Gamma$-periodic solution to the equation \eqref{1.8}, then $\tilde\lambda=\lambda+\lambda_0$ is a positive function for a sufficiently large constant $\lambda_0$ and  $\tilde\lambda$ is again a $\Gamma$-periodic solution to \eqref{1.8}.
However, the answer to the question is negative. In author's opinion, the following theorem is the main result of the present work.

\begin{theorem} \label{Th1.2}
The following statement is true for an arbitrary lattice $\Gamma\subset\R^2$. Let $\lambda\in C^\infty(\R^2)$ be a real $\Gamma$-periodic solution to the system \eqref{1.6}. Assume that $\lambda+\lambda_0$ is also a solution to the system \eqref{1.6} with some real constant $\lambda_0\neq0$. Then $\lambda$ is a one-dimensional solution.
\end{theorem}

We explain the scheme of the proof of Theorem \ref{Th1.2}. The following statement is just an easy remark.

\begin{proposition} \label{P1.1}
Let $\lambda\in C^3(\R^2)$ be a real solution to the equation \eqref{1.8}.  Assume that $\lambda+\lambda_0$ is also a solution to the equation \eqref{1.8} with some real constant $\lambda_0\neq0$. Then $\lambda$ solves also the equation
\begin{equation}
c\lambda_{zzz}+\bar c\lambda_{\bar z\bar z\bar z}=0.
                         \label{1.12}
\end{equation}
Conversely, if a real function $\lambda\in C^3(\R^2)$ solves equations \eqref{1.8} and \eqref{1.12}, then $\lambda+\lambda_0$ also solves the equation \eqref{1.8} for any real constant $\lambda_0$.
\end{proposition}

The linear equation \eqref{1.12} can be easily investigated. In particular, the following statement holds.

\begin{proposition} \label{P1.2}
For every lattice $\Gamma\subset\R^2$ and for every complex constant $c\neq0$ the following statement is true.
A $\Gamma$-periodic function $\lambda\in C^3(\R^2)$ solves the equation \eqref{1.12} if and only if its spectrum is contained in the union of three straight lines passing through the origin and intersecting with each other at equal angles.
\end{proposition}

The system \eqref{1.6} is essentially simplified with the help of the latter proposition. The two-dimensional array $\{{\hat\lambda}_k\}_{k\in\Gamma'}$ of unknowns is replaced with three sequences $\{x_n\}_{n\in\Z}$, $\{y_n\}_{n\in\Z}$, $\{z_n\}_{n\in\Z}$ corresponding to three lines in Proposition \ref{P1.2}. Besides this, five above-mentioned constants can be eliminated from the system. In this way Theorem \ref{Th1.2} is reduced to the following statement.

\begin{theorem} \label{Th1.3}
Let three sequences of complex numbers $\{x_n\}_{n\in\Z\setminus\{0\}}$, $\{y_n\}_{n\in\Z\setminus\{0\}}$,\\ $\{z_n\}_{n\in\Z\setminus\{0\}}$ constitute a solution to the system
\begin{equation}
x_{n_1}y_{n_2}-x_{n_1+n_2}z_{n_2}+y_{n_1+n_2}z_{-n_1}=0\quad(n_1\neq0,n_2\neq0,n_1+n_2\neq0).
                         \label{1.13}
\end{equation}
Assume additionally the solution to satisfy the parity condition
\begin{equation}
x_{-n}=\overline{x_n},\quad y_{-n}=\overline{y_n},\quad z_{-n}=\overline{z_n}
                         \label{1.14}
\end{equation}
and the decay condition
\begin{equation}
\sum\limits_{n=1}^\infty(|x_{n}|+|y_{n}|+|z_{n}|)<\infty.
                         \label{1.15}
\end{equation}
Then it is a one-dimensional solution, i.e., one of the following three statements holds:
{\rm (1)} $y_n=z_n=0$ for all $n\neq0$,
{\rm (2)} $x_n=z_n=0$ for all $n\neq0$,
{\rm (3)} $x_n=y_n=0$ for all $n\neq0$.
\end{theorem}

{\bf Remark.} The author does not know whether Theorem \ref{Th1.3} remains valid if the hypothesis \eqref{1.15} is either deleted or replaced with the weaker one:
$|x_{n}|+|y_{n}|+|z_{n}|\rightarrow0$ as $n\rightarrow\infty$.

In Section 2, we prove Theorem \ref{Th1.1} and derive the equation \eqref{1.8} from the latter theorem. We also present an alternative version of the system  \eqref{1.6}, where Fourier coefficients are numbered by pairs of integers but not by nodes of a lattice. Theorem \ref{Th1.2} is reduced to Theorem \ref{Th1.3} in Section 3. Last two sections are devoted to the proof of Theorem \ref{Th1.3} which turns out to be quite difficult.

\section{Proof of Theorem \ref{Th1.1}}

As a start point, we use the following statement proved in \cite[Lemma 4.1]{Sh}.

\begin{proposition} \label{P2.1}
A Riemannian torus $({\mathbb T}^2,g)=\big({\mathbb R}^2/\Gamma,\lambda(dx^2+dy^2)\big)$ admits a non-trivial third rank Killing tensor field if and only if, for some pair $c=(c_1,c_2)\neq(0,0)$ of real constants, the equation
\begin{equation}
{\nabla}{\nabla}u-\frac{1}{2}(\Delta_g u)g=T^c
                         \label{2.1}
\end{equation}
has a solution $u\in C^\infty({\mathbb R}^2)$ with $\Gamma$-periodic derivatives
$u_x$ and $u_y$. Here $\Delta_g$ is the Laplace -- Beltrami operator of the metric $g$ and $T^c$ is the symmetric second rank tensor field on the torus whose coordinates are defined by
$$
T^c_{11}=-T^c_{22}=\frac{1}{2}\lambda(-c_2\lambda_x+c_1\lambda_y),\quad
T^c_{12}=\frac{1}{2}\lambda(c_1\lambda_x+c_2\lambda_y).
$$
\end{proposition}

(Unfortunately, there is the inaccuracy in the statement of Lemma 4.1 in \cite{Sh}: the word ``irreducible'' is written instead of ``non-trivial''. I am sure the reader will easily understand this simple statement.)

Computing the coordinates of the tensor field from the left-hand side of \eqref{2.1} with the help of standard tensor analysis formulas, we make sure that this equation is equivalent to the system
$$
\begin{aligned}
u_{xx}-u_{yy}-\lambda^{-1}\lambda_xu_x+\lambda^{-1}\lambda_yu_y&=\lambda(-c_2\lambda_x+c_1\lambda_y),\\
2u_{xy}-\lambda^{-1}\lambda_yu_x-\lambda^{-1}\lambda_xu_y&=\lambda(c_1\lambda_x+c_2\lambda_y).
\end{aligned}
$$
We are interested in the consistency conditions for the system. It can be written in the form
$$
\frac{\partial(\lambda^{-1}u_x)}{\partial x}-\frac{\partial(\lambda^{-1}u_y)}{\partial y}=-c_2\lambda_x+c_1\lambda_y,\quad
\frac{\partial(\lambda^{-1}u_x)}{\partial y}+\frac{\partial(\lambda^{-1}u_y)}{\partial x}=c_1\lambda_x+c_2\lambda_y.
$$
Introducing the functions
\begin{equation}
v=\lambda^{-1}u_x,\quad w=\lambda^{-1}u_y,
                         \label{2.3}
\end{equation}
we arrive to the Cauchy -- Riemann equations
\begin{equation}
v_x-w_y=-c_2\lambda_x+c_1\lambda_y,\quad v_y+w_x=c_1\lambda_x+c_2\lambda_y.
                         \label{2.4}
\end{equation}
By \eqref{2.3}, the consistency condition holds
\begin{equation}
\frac{\partial(\lambda v)}{\partial y}-\frac{\partial(\lambda w)}{\partial x}=0.
                         \label{2.5}
\end{equation}

Since $u_x,u_y$ and $\lambda$ are $\Gamma$-periodic functions, $v$ and $w$ are also $\Gamma$-periodic functions. Recall that we have written down the Fourier series for $\lambda$, see \eqref{1.5}. Let us write similar formulas for $v$ and $w$:
\begin{equation}
v=\sum\limits_{n\in\Gamma'}{\hat v}_n\,e^{i(n_1x+n_2y)},\quad w=\sum\limits_{n\in\Gamma'}{\hat w}_n\,e^{i(n_1x+n_2y)}.
                         \label{2.6}
\end{equation}
Inserting these expressions into equations \eqref{2.4} and performing termwise differentiation of Fourier series, we obtain
$$
\begin{aligned}
\sum\limits_{n\in\Gamma'}(n_1{\hat v}_n-n_2{\hat w}_n)e^{i(n_1x+n_2y)}&=\sum\limits_{n\in\Gamma'}(-c_2n_1+c_1n_2){\hat\lambda}_ne^{i(n_1x+n_2y)},\\
\sum\limits_{n\in\Gamma'}(n_2{\hat v}_n+n_1{\hat w}_n)e^{i(n_1x+n_2y)}&=\sum\limits_{n\in\Gamma'}(c_1n_1+c_2n_2){\hat\lambda}_ne^{i(n_1x+n_2y)}.
\end{aligned}
$$
Equating coefficients at the same exponents, we arrive to the system
$$
n_1{\hat v}_n-n_2{\hat w}_n=(-c_2n_1+c_1n_2){\hat\lambda}_n,\quad
n_2{\hat v}_n+n_1{\hat w}_n=(c_1n_1+c_2n_2){\hat\lambda}_n\quad(n\in\Gamma').
$$
The system is uniquely solvable for $n\neq0$:
$$
{\hat v}_n=\frac{-c_2n_1^2+2c_1n_1n_2+c_2n_2^2}{|n|^2}\,{\hat\lambda}_n,\quad
{\hat w}_n=\frac{c_1n_1^2+2c_2n_1n_2-c_1n_2^2}{|n|^2}\,{\hat\lambda}_n.
$$
The Fourier series \eqref{2.6} take now the form
$$
\begin{aligned}
v&={\hat v}_0+\sum\limits_{k\in\Gamma'\setminus\{0\}}\frac{-c_2k_1^2+2c_1k_1k_2+c_2k_2^2}{|k|^2}\,{\hat\lambda}_k\,e^{i(k_1x+k_2y)},\\
w&={\hat w}_0+\sum\limits_{k\in\Gamma'\setminus\{0\}}\frac{c_1k_1^2+2c_2k_1k_2-c_1k_2^2}{|k|^2}\,{\hat\lambda}_k\,e^{i(k_1x+k_2y)}.
\end{aligned}
$$
Using these formulas and \eqref{1.5}, we find Fourier series for the functions $\lambda v$ and $\lambda w$:
$$
\begin{aligned}
\lambda v&=\sum\limits_{n\in\Gamma'}\Big({\hat v}_0{\hat\lambda}_n
+\sum\limits_{k\in\Gamma'\setminus\{0\}}\frac{-c_2k_1^2+2c_1k_1k_2+c_2k_2^2}{|k|^2}\,{\hat\lambda}_k{\hat\lambda}_{n-k}\Big)e^{i(n_1x+n_2y)},\\
\lambda w&=\sum\limits_{n\in\Gamma'}\Big({\hat w}_0{\hat\lambda}_n
+\sum\limits_{k\in\Gamma'\setminus\{0\}}\frac{c_1k_1^2+2c_2k_1k_2-c_1k_2^2}{|k|^2}\,{\hat\lambda}_k{\hat\lambda}_{n-k}\Big)e^{i(n_1x+n_2y)}.
\end{aligned}
$$
Substitute these expressions into the equation \eqref{2.5}
$$
\begin{aligned}
\frac{\partial}{\partial y}&\sum\limits_{n\in\Gamma'}\Big({\hat v}_0{\hat\lambda}_n
+\sum\limits_{k\in\Gamma'\setminus\{0\}}\frac{-c_2k_1^2+2c_1k_1k_2+c_2k_2^2}{|k|^2}\,{\hat\lambda}_k{\hat\lambda}_{n-k}\Big)e^{i(n_1x+n_2y)}\\
-\frac{\partial}{\partial x}&\sum\limits_{n\in\Gamma'}\Big({\hat w}_0{\hat\lambda}_n
+\sum\limits_{k\in\Gamma'\setminus\{0\}}\frac{c_1k_1^2+2c_2k_1k_2-c_1k_2^2}{|k|^2}\,{\hat\lambda}_k{\hat\lambda}_{n-k}\Big)e^{i(n_1x+n_2y)}=0.
\end{aligned}
$$
After termwise differentiation of series, the equation becomes
$$
\begin{aligned}
&\sum\limits_{n\in\Gamma'}\Big(({\hat v}_0n_2-{\hat w}_0n_1){\hat\lambda}_n\\
+&\sum\limits_{k\in\Gamma'\setminus\{0\}}\frac{n_2(-c_2k_1^2+2c_1k_1k_2+c_2k_2^2)-n_1(c_1k_1^2+2c_2k_1k_2-c_1k_2^2)}{|k|^2}\,
{\hat\lambda}_k{\hat\lambda}_{n-k}\Big)e^{i(n_1x+n_2y)}=0.
\end{aligned}
$$
Equating coefficients of the series on the left-hand side to zero, we arrive to equations \eqref{1.6}, where $a_1={\hat w}_0,a_2=-{\hat v}_0$. It remains to notice that $a_1$ and $a_2$ are real constants since $v$ and $w$ are real functions. Theorem \ref{Th1.1} is proved.

\bigskip

The equivalence of the equation \eqref{1.8} and system \eqref{1.6} is proved by a straightforward calculation. To this end one finds the Fourier series of the function from the left-hand side of \eqref{1.8} on using \eqref{1.5}, \eqref{1.7} and well known formulas
$$
\partial_z=\frac{1}{2}(\partial_x-i\partial_y),\quad\partial_{\bar z}=\frac{1}{2}(\partial_x+i\partial_y).
$$
Then, equating coefficients of the resulting Fourier series to zero, we again arrive to \eqref{1.6}. These arguments are invertible.

In formulas \eqref{1.5} and \eqref{1.8}, Fourier coefficients of the function $\lambda$ are numbered by nodes of the lattice $\Gamma'$. How do these formulas look like if the Fourier coefficients are numbered by pairs of integers? Let $(e_1,e_2)$ be a basis of the lattice $\Gamma'$. Since our problem is invariant under rotations and homotheties of the lattice, we can assume without lost of generality that $e_1=(1,0)$. Let $e_2=(b,d)$. We can assume that $d>0$ since $(e_1,-e_2)$ is also a basis. The pair $(e_1,ke_1+e_2)$ is again a basis of $\Gamma'$ for every integer $k$ and $ke_1+e_2=(b+k,d)$. Therefore we can assume without lost of generality that $b\in[0,1)$. Thus, the lattice $\Gamma'$ is determined by two parameters $b\in[0,1)$ and $d>0$. The vectors $e_1=(1,0)$ and $e_2=(b,d)$ constitute the basis of $\Gamma'$. Every node $m=(m_1,m_2)\in\Gamma'$ is uniquely represented in the form
\begin{equation}
m=n_1e_1+n_2e_2=(n_1+bn_2,dn_2)\quad(n_1,n_2\in\Z).
                         \label{2.7}
\end{equation}
Let us reproduce formula \eqref{1.5}
$$
\lambda(x,y)=\sum\limits_{m\in\Gamma'}{\hat\lambda}_{m_1,m_2}e^{i(m_1x+m_2y)}.
$$
After the change \eqref{2.7}, the formula takes the form
$$
\lambda(x,y)=\sum\limits_{n\in\Z^2}{\hat\lambda}_{n_1+bn_2,dn_2}e^{i((n_1+bn_2)x+dn_2y)}.
$$
Let us introduce the notation that will be used in the rest of the paper
$$
{\check\lambda}_{n_1,n_2}={\hat\lambda}_{n_1+bn_2,dn_2}\quad(n_1,n_2\in\Z).
$$
Thus, ${\check\lambda}_{n_1,n_2}$ are Fourier coefficients of the function $\lambda$ numbered by pairs of integers. The Fourier series of a $\Gamma$-periodic function $\lambda$ looks now as follows:
$$
\lambda(x,y)=\sum\limits_{(n_1,n_2)\in\Z^2}{\check\lambda}_{n_1,n_2}e^{i((n_1+bn_2)x+dn_2y)}.
$$

By the same change \eqref{2.7}, the system \eqref{1.8} is transformed to the form
\begin{equation}
\sum\limits_{(k_1,k_2)\in{\Z}^2\setminus\{(0,0)\}}
\psi(n_1,n_2;k_1,k_2;c_1,c_2;b,d)\,{\check\lambda}_{k_1,k_2}{\check\lambda}_{n_1-k_1,n_2-k_2}
=(a_1n_1+a_2n_2){\check\lambda}_{n_1,n_2},
                         \label{2.9}
\end{equation}
where
\begin{equation}
\psi(n_1,n_2;k_1,k_2;c_1,c_2;b,d)=c_1\psi_1(n_1,n_2;k_1,k_2;b,d)+c_2\psi_2(n_1,n_2;k_1,k_2;b,d),
                         \label{2.10}
\end{equation}
\begin{equation}
\psi_1(n_1,n_2;k_1,k_2;b,d)=d\,\frac{n_2k_1^2+2(n_1+2 b n_2)k_1k_2+\big(2 b n_1+(3 b^2- d^2)n_2\big)k_2^2}{(k_1+ b k_2)^2+ d^2k_2^2},
                         \label{2.11}
\end{equation}
\begin{equation}
\begin{aligned}
&\psi_2(n_1,n_2;k_1,k_2;b,d)=\\
&=\frac{-(n_1\!+\! b n_2)k_1^2+2\big(- b n_1\!+\!( d^2\!-\! b^2)n_2\big)k_1k_2
+\big(( d^2\!-\! b^2)n_1+ b(3 d^2\!-\! b^2)n_2\big)k_2^2}{(k_1+ b k_2)^2+ d^2k_2^2}.
\end{aligned}
                         \label{2.12}
\end{equation}
At the same time, we had to change the value of the constant $a_2$ in accordance with the formula $a_2:=ba_1+da_2$. However, the pair of constants $(a_1,ba_1+da_2)$ is as arbitrary as the pair $(a_1,a_2)$.

\section{Proof of Theorem \ref{Th1.2}}

We first present the proof of Proposition \ref{P1.1}.

Let $\lambda\in C^3(\R^2)$ be a real solution to the equation \eqref{1.8}. Assume that the function $\lambda+\lambda_0$ also solves the equation \eqref{1.8} with some real constant $\lambda_0\neq0$, i.e.,
$$
\frac{\partial}{\partial z}\Big((\lambda+\lambda_0)\big(c\Delta^{-1}(\lambda+\lambda_0)_{zz}+a\big)\Big)
+\frac{\partial}{\partial\bar z}\Big((\lambda+\lambda_0)\big(\bar c\Delta^{-1}(\lambda+\Lambda_0)_{\bar z\bar z}+\bar a\big)\Big)=0.
$$
After the obvious simplification, this takes the form
$$
\frac{\partial}{\partial z}\Big(\lambda(c\Delta^{-1}\lambda_{zz}+a)\Big)
+\frac{\partial}{\partial\bar z}\Big(\lambda(\bar c\Delta^{-1}\lambda_{\bar z\bar z}+\bar a)\Big)
+\lambda_0\Delta^{-1}(c\lambda_{zzz}+\bar c\lambda_{\bar z\bar z\bar z})=0.
$$
Subtracting the equation \eqref{1.8} from this equality, we obtain $\Delta^{-1}(c\lambda_{zzz}+\bar c\lambda_{\bar z\bar z\bar z})=0$. Applying the operator $\Delta$ to the latter equation, we arrive to \eqref{1.12}. These arguments are invertible.

Next we present the proof of Proposition \ref{P1.2}.

Differentiating the Fourier series \eqref{1.5}, we have
$$
\lambda_{zzz}=-\frac{i}{8}\sum\limits_{n\in\Gamma'}(n_1-in_2)^3\hat\lambda_n\,e^{i(n_1x+n_2y)},\quad
\lambda_{\bar z\bar z\bar z}=-\frac{i}{8}\sum\limits_{n\in\Gamma'}(n_1+in_2)^3\hat\lambda_n\,e^{i(n_1x+n_2y)}.
$$
Substituting this expression into \eqref{1.12}, we arrive to the equation
$$
\Re\big(c(n_1-in_2)^3\big)=0
$$
which should hold for nodes $(n_1,n_2)$ of the lattice $\Gamma'$ belonging to the spectrum of the function $\lambda$.
If $c=c_1+ic_2$, then the latter equation can be written in the form
\begin{equation}
c_1(n_1^2-3n_2^2)n_1+c_2(3n_1^2-n_2^2)n_2=0.
                         \label{3.1}
\end{equation}
Without lost of generality we can assume that $c_1=-\sin3\alpha,\ c_2=\cos3\alpha$ in the equation \eqref{3.1} with some real $\alpha$. We set also $n_1=\sqrt{n_1^2+n_2^2}\cos\varphi,\ n_2=\sqrt{n_1^2+n_2^2}\sin\varphi$. Then the equation \eqref{3.1} takes the form
$\sin(3\varphi-3\alpha)=0$. Its solutions are: $\varphi=k\frac{\pi}{3}+\alpha\ (k\in\Z)$. This proves Proposition \ref{P1.2}.

\bigskip

We proceed to the proof of Theorem \ref{Th1.2} assuming Theorem \ref{Th1.3} to be valid.

Let a lattice $\Gamma\subset\R^2$ and real $\Gamma$-periodic function $\lambda\in C^\infty(\R^2)$ satisfy hypotheses of Theorem \ref{Th1.2}. By Proposition \ref{P1.2}, the spectrum of the function $\lambda$ is contained in the union of three lines $L_0\cup L_1\cup L_2$ passing through the origin and intersecting with each other at the angle $\pi/3$. We have to prove that $\lambda$ is a one-dimensional solution, i.e., the spectrum of $\lambda$ is contained in one of these lines. Contrary to the assertion of the theorem, suppose that each of the lines contains at least one point belonging to the spectrum of $\lambda$ and different of the origin (the case when the spectrum lies on two lines is considered in a similar way with many simplifications). Using the invariance of our problem under rotations, we can assume that
\begin{equation}
L_0=\{(x,0)\mid x\in\R\},\quad L_1=\{(x,\sqrt{3}x)\mid x\in\R\},\quad L_2=\{(x,-\sqrt{3}x)\mid x\in\R\}.
                         \label{3.2}
\end{equation}
Since the line $L_0$ contains a point of the spectrum of $\lambda$ other than the origin, there are non-zero nodes of the lattice $\Gamma'$ on the line $L_0$. Let $e_1=(x_0,0)\ (x_0>0)$  be the closest to the origin point of the kind $(x,0)\in\Gamma'\cap L_0\ (x>0)$. Using the invariance of our problem with respect to homotheties centered at the origin, we can assume that $x_0=1$, i.e., $e_1=(1,0)$. Let us demonstrate that the lattice $\Gamma'$ has a basis of the form
\begin{equation}
e_1=(1,0),\quad e_2=(b,d)\quad\big(b\in[0,1), d>0\big).
                         \label{3.3}
\end{equation}
Indeed, let $\big(e'_1=(e'^1_1,e'^2_1),e'_2=(e'^1_2,e'^2_2)\big)$ be an arbitrary basis of $\Gamma'$. We represent the vector $e_1$ as $e_1=n_1e'_1+n_2e'_2\ (n_1,n_2\in\Z)$. Then $(n_1,n_2)\neq0$ and $(n_1,n_2)$ are relatively prime integers since $e_1$ is the point of the set $L_0\cap\Gamma'$ closest to the origin, other than the origin. There exists a pair of integers $(l_1,l_2)$ such that $n_1l_2-n_2l_1=1$. The matrix $A=\left(\begin{array}{cc}n_1&n_2\\l_1&l_2\end{array}\right)$ belongs to $SL(2,\Z)$ and hence the pair of vectors
$$
{{e_1}\choose{e''_2}}=A{{e'_1}\choose{e'_2}}
$$
is a basis of the lattice $\Gamma'$. Starting with the basis $(e_1,e''_2)$ and repeating our arguments presented before formula \eqref{2.7}, we find a basis of $\Gamma'$ of the form \eqref{3.3}.

Let $(x,\sqrt{3}x)$ be the node of the lattice $\Gamma'$ closest to the origin, other then the origin and such that $x>0$. This means the existence of the pair $(p_1,p_2)$ of relatively prime integers such that $(x,\sqrt{3}x)=p_1(1,0)+p_2(b,d)$, i.e.,
\begin{equation}
x=p_1+p_2b,\quad \sqrt{3}x=p_2d\quad(p_2>0).
                         \label{3.4}
\end{equation}
Similarly, let $(x',-\sqrt{3}x')$ be the node of the lattice $\Gamma'$ closest to the origin, other then the origin and such that $x'>0$. This means the existence of the pair $(q_1,q_2)$ of relatively prime integers such that $(x',-\sqrt{3}x')=q_1(1,0)+q_2(b,d)$, i.e.,
\begin{equation}
x'=q_1+q_2b,\quad -\sqrt{3}x'=q_2d\quad(q_2>0).
                         \label{3.5}
\end{equation}

From \eqref{3.4}--\eqref{3.5}
$$
-b+d/\sqrt{3}=p_1/p_2,\quad -b-d/\sqrt{3}=q_1/q_2.
$$
Solving the system, we have
\begin{equation}
b=-\frac{1}{2}\Big(\frac{p_1}{p_2}+\frac{q_1}{q_2}\Big),\quad d=\frac{\sqrt{3}}{2}\Big(\frac{p_1}{p_2}-\frac{q_1}{q_2}\Big).
                         \label{3.7}
\end{equation}

Equalities \eqref{3.7} allow us to eliminate the constants $b$ and $d$ from all our formulas. The lattice $\Gamma'$ is now determined by two pairs $(p_1,p_2>0)$ and $(q_1,q_2>0)$ of relatively prime integers. But the pairs should satisfy some inequalities. Indeed, the conditions $d>0$ and $b\in[0,1)$ mean that
\begin{equation}
p_1q_2-p_2q_1>0,\quad -2p_2q_2<p_1q_2+p_2q_1\leq0.
                         \label{3.8}
\end{equation}
This is equivalent to the system of inequalities
\begin{equation}
p_2>0,\quad q_2>0,\quad q_1<0,\quad\max\Big\{\frac{q_1}{q_2}\,p_2,-\big(2+\frac{q_1}{q_2}\big)\,p_2\Big\}<p_1\leq-\frac{q_1}{q_2}\,p_2.
                         \label{3.9}
\end{equation}
One can easily prove the converse statement: If two pairs of relatively prime integers $(p_1,p_2)$ and $(q_1,q_2)$ satisfy \eqref{3.9}, then the reals $(b,d)$, defined by \eqref{3.7}, satisfy $b\in[0,1),\ d>0$ and the lattice $\Gamma'=\{n_1(1,0)+n_2(b,d)\mid(n_1,n_2)\in\Z^2\}$ has the following intersections with the lines \eqref{3.2}:
$$
\begin{aligned}
L_0\cap\Gamma'&=\{(n,0)\mid n\in\Z\},\quad
L_1\cap\Gamma'=\{\big(n(p_1+bp_2),ndp_2\big)\mid n\in\Z\},\\
L_2\cap\Gamma'&=\{\big(n(q_1+bq_2),ndq_2\big)\mid n\in\Z\}.
\end{aligned}
$$
In the case of $(p_1,p_2)=(0,1),\ (q_1,q_2)=(-1,1)$, the lattice $\Gamma'$ is the result of tiling the plane with regular triangles. This grid is sometimes called a {\it honeycomb}. The author used this lattice to control many of formulas below. Although the honeycomb is not mentioned in the rest of the paper, the reader is recommended to keep this simplest example in his/her mind.

Since the spectrum of the function $\lambda$ is contained in the union of three lines \eqref{3.2}, the Fourier coefficients ${\check\lambda}_{n_1,n_2}$, introduced in the previous section, are expressed by
\begin{equation}
{\check\lambda}_{n_1,n_2}=\left\{\begin{array}{l}
\alpha_n\ \mbox{if}\ (n_1,n_2)=(n,0),\\
\beta_n\ \mbox{if}\ (n_1,n_2)=(np_1,np_2),\\
\gamma_n\ \mbox{if}\ (n_1,n_2)=(nq_1,nq_2),\\
0\ \mbox{otherwise},\end{array}\right.
                         \label{3.10}
\end{equation}
where $\{\alpha_n\}_{n\in\Z},\{\beta_n\}_{n\in\Z},\{\gamma_n\}_{n\in\Z}$ are three sequences of complex numbers which are considered as unknowns in what follows. In terms of these unknowns, the parity condition \eqref{1.10} is written as
\begin{equation}
\alpha_{-n}=\overline{\alpha_n},\quad
\beta_{-n}=\overline{\beta_n},\quad
\gamma_{-n}=\overline{\gamma_n};
                         \label{3.11}
\end{equation}
and the decay condition \eqref{1.11} takes the form
\begin{equation}
|\alpha_{n}|+|\beta_{n}|+|\gamma_n|\le C_N(|n|+1)^{-N}\quad\mbox{for every}\ n\in\N.
                         \label{3.12}
\end{equation}
Besides this, the following equalities should hold:
$$
\alpha_0=\beta_0=\gamma_0={\check\lambda}_{0,0}.
$$

Due to \eqref{3.10}, the system \eqref{2.9} admits several simplifications. First of all the double sum is replaced with three single sums, i.e., the equation \eqref{2.9} takes the form
\begin{equation}
\begin{aligned}
&\sum\limits_{k\in\Z\setminus\{0\}}
\psi(n_1,n_2;k,0;c_11,c_2;b,d)\,\alpha_k{\check\lambda}_{n_1-k,n_2}\\
+&\sum\limits_{k\in\Z\setminus\{0\}}
\psi(n_1,n_2;p_1k,p_2k;c_1,c_2;b,d)\,\beta_k{\check\lambda}_{n_1-p_1k,n_2-p_2k}\\
+&\sum\limits_{k\in\Z\setminus\{0\}}
\psi(n_1,n_2;q_1k,q_2k;c_1,c_2;b,d)\,\gamma_k{\check\lambda}_{n_1-q_1k,n_2-q_2k}
=(a_1n_1+a_2n_2){\check\lambda}_{n_1,n_2}.
\end{aligned}
                         \label{3.14}
\end{equation}
Coefficients of the equation \eqref{3.14} are expressed by
\begin{equation}
\psi(n_1,n_2;k,0;c_1,c_2;b,d)=-2c^2n_1+\Big(c_1+c_2\big(\frac{p_1}{p_2}+\frac{q_1}{q_2}\big)\Big)n_2,
                         \label{3.15}
\end{equation}
\begin{equation}
\begin{aligned}
&\varphi(n_1,n_2;p_1,p_2;q_1,q_2;c_1,c_2):=\psi(n_1,n_2;p_1k,p_2k;c_1,c_2;b,d)\\
&=c_1q_2\,\frac{p_2n_1-p_1n_2}{p_1q_2-p_2q_1}+c_2\Big(n_1+
\frac{\frac{p_1^3}{p_2}-4\frac{q_1}{q_2}p_1^2+5\frac{q_1^2}{q_2^2}p_1p_2-2\frac{q_1^3}{q_2^3}p_2^2}{\big(p_1-\frac{q_1}{q_2}p_2\big)^2}\,n_2\Big),
\end{aligned}
                         \label{3.16}
\end{equation}
\begin{equation}
\begin{aligned}
&\varphi(n_1,n_2;q_1,q_2;p_1,p_2;c_1,c_2):=\psi(n_1,n_2;q_1k,q_2k;c_1,c_2;b,d)\\
&=-c_1p_2\,\frac{q_2n_1-q_1n_2}{p_1q_2-p_2q_1}+c_2\Big(n_1+\frac{
\frac{q_1^3}{q_2}-4\frac{p_1}{p_2}q_1^2+5\frac{p_1^2}{p_2^2}q_1q_2-2\frac{p_1^3}{p_2^3}q_2^2}{\big(q_1-\frac{p_1}{p_2}q_2\big)^2}\,n_2\Big).
\end{aligned}
                         \label{3.17}
\end{equation}
These formulas are obtained from \eqref{2.10}--\eqref{2.12} by substituting values \eqref{3.7} for the parameters $(b,d)$ and by substituting the values $(k_1,k_2)=(k,0)$, $(k_1,k_2)=(p_1k,p_2k)$ and $(k_1,k_2)=(q_1k,q_2k)$ respectively. At the same time we have changed the values of constants as $c_1:=c_1/d,c_2:=2c_2$. The most important (although quite obvious) circumstance is that the coefficients \eqref{3.15}--\eqref{3.17} are independent of the summation variable $k$. In virtue of the circumstance, the equation \eqref{3.14} can be rewritten as
\begin{equation}
\begin{aligned}
&\Big(-2c_2n_1+\big(c_1+c_2\big(\frac{p_1}{p_2}+\frac{q_1}{q_2}\big)\big)n_2\Big)\sum\limits_{k\in\Z\setminus\{0\}}\alpha_k{\check\lambda}_{n_1-k,n_2}\\
+&\varphi(n_1,n_2;p_1,p_2;q_1,q_2;c_1,c_2)\sum\limits_{k\in\Z\setminus\{0\}}\beta_k{\check\lambda}_{n_1-p_1k,n_2-p_2k}\\
+&\varphi(n_1,n_2;q_1,q_2;p_1,p_2;c_1,c_2)\sum\limits_{k\in\Z\setminus\{0\}}\gamma_k{\check\lambda}_{n_1-q_1k,n_2-q_2k}
=(a_1n_1+a_2n_2){\check\lambda}_{n_1,n_2}.
\end{aligned}
                         \label{3.18}
\end{equation}

Possible simplification are still not exhausted. The second factor of each summand on the left-hand side of \eqref{3.18} can be expressed through $(\alpha,\beta,\gamma)$, as well as the right-hand side of the equation. To this end we have to consider separately four possible cases corresponding four lines on the right-hand side of the formula \eqref{3.10}.

First of all we observe that the equation \eqref{3.18} tautologically holds in the case of $(n_1,n_2)=(0,0)$ (both sides of the equation are equal to zero). This case is excluded from our further considerations.

We first set $(n_1,n_2)=(n,0)$ in \eqref{3.18}, where $n\neq0$. In other words, we consider the equation \eqref{3.18} when the node of the lattice $\Gamma'$, numbered by the pair $(n_1,n_2)$, belongs to the line $L_0$. In this case the equation \eqref{3.18} is of the form
\begin{equation}
\begin{aligned}
&-2nc_2\sum\limits_{k\in\Z\setminus\{0\}}\alpha_k\alpha_{n-k}\\
&+\varphi(n,0;p_1,p_2;q_1,q_2;c_1,c_2)\sum\limits_{k\in\Z\setminus\{0\}}\beta_k{\check\lambda}_{n-p_1k,-p_2k}\\
&+\varphi(n,0;q_1,q_2;p_1,p_2;c_1,c_2)\sum\limits_{k\in\Z\setminus\{0\}}\gamma_k{\check\lambda}_{n-q_1k,-q_2k}
=na_1\alpha_n.
\end{aligned}
                         \label{3.19}
\end{equation}

Let us pay our attention to the zeroth Fourier coefficient ${\check\lambda}_{0,0}=\alpha_0=\beta_0=\gamma_0$. By the second statement of Proposition \ref{P1.1}, the coefficient cannot participate in the system \eqref{3.18} since it plays the role of the constant $\lambda_0$ from Proposition \ref{P1.1}. On the other hand, in the case of $c_2\neq0$, there is the term $-2nc_2\alpha_n\alpha_0$ on the left-hand side of the equation \eqref{3.19} which does not cancel with other terms. It remains to observe that all $\alpha_n\ (n\neq0)$ cannot be equal to zero since the line $L_0$ contains at least one point from the spectrum of $\lambda$ which is not the origin. We thus arrive to the important conclusion: $c_2=0$. Since $(c_1,c_2)\neq(0,0)$ and only pairwise ratios $(a_1:a_2:c_1:c_2)$ are essential, we can assume without lost of generality that $(c_1,c_2)=(1,0)$.

Formulas \eqref{3.16}--\eqref{3.17} are now simplified to the following one:
\begin{equation}
\varphi(n_1,n_2;p_1,p_2;q_1,q_2)=q_2\,\frac{p_2n_1-p_1n_2}{p_1q_2-p_2q_1}.
                         \label{3.20}
\end{equation}
Of course the constants $(c_1,c_2)=(1,0)$ do not participate in the list of arguments anymore. In particular,
$$
\varphi(n,0;p_1,p_2;q_1,q_2)=-\varphi(n,0;q_1,q_2;p_1,p_2)=\frac{p_2q_2}{p_1q_2-p_2q_1}\,n
$$
and the equation \eqref{3.19} takes the form
\begin{equation}
\sum\limits_{k\in\Z\setminus\{0\}}\beta_k{\check\lambda}_{n-p_1k,-p_2k}
-\sum\limits_{k\in\Z\setminus\{0\}}\gamma_k{\check\lambda}_{n-q_1k,-q_2k}
=a_1\frac{p_1q_2-p_2q_1}{p_2q_2}\,\alpha_n\quad(n\neq0).
                         \label{3.21}
\end{equation}

We first analyze the first sum on the left-hand side of \eqref{3.21}. Its summand $\beta_k{\check\lambda}_{n-p_1k,-p_2k}$ can be nonzero in three cases only:

(1) ${\check\lambda}_{n-p_1k,-p_2k}=\alpha_m$ if $(n-p_1k,-p_2k)=(m,0)$ for some integer $m$,

(2) ${\check\lambda}_{n-p_1k,-p_2k}=\beta_m$ if $(n-p_1k,-p_2k)=(mp_1,mp_2)$ for some integer $m$,

(3) ${\check\lambda}_{n-p_1k,-p_2k}=\gamma_m$ if $(n-p_1k,-p_2k)=(mq_1,mq_2)$ for some integer $m$.

The first case is impossible since $p_2k\neq0$. Let us demonstrate that the second case is also impossible. Indeed, in such a case
$$
n-p_1k=mp_1,\quad -p_2k=mp_2.
$$
From this $n=0$ that contradicts the assumption $n\neq0$. Thus, the third case remains when
$$
n-p_1k=mq_1,\quad -p_2k=mq_2.
$$
From this
$$
k=\frac{q_2}{p_1q_2-p_2q_1}\,n,\quad m=-\frac{p_2}{p_1q_2-p_2q_1}\,n.
$$
Therefore the integers $p_2n$ and $q_2n$ should be divisible by $p_1q_2-p_2q_1$. Thus,
\begin{equation}
\sum\limits_{k\in\Z\setminus\{0\}}\beta_k{\check\lambda}_{n-p_1k,-p_2k}=
\left\{\begin{array}{l}\beta_r\gamma_s\ \mbox{if}\ q_2n=r(p_1q_2-p_2q_1)\\
\qquad \mbox{and}\ p_2n=-s(p_1q_2-p_2q_1)\ (r,s\in\Z);\\
0\ \mbox{otherwise}.\end{array}\right.
                         \label{3.22}
\end{equation}

Next, we analyze the second sum on the left-hand side of \eqref{3.21}. Its summand $\gamma_k{\check\lambda}_{n-q_1k,-q_2k}$ can be nonzero in three cases only:

(1) ${\check\lambda}_{n-q_1k,-q_2k}=\alpha_m$ if $(n-q_1k,-q_2k)=(m,0)$ for some integer $m$,

(2) ${\check\lambda}_{n-q_1k,-q_2k}=\beta_m$ if $(n-q_1k,-q_2k)=(mp_1,mp_2)$ for some integer $m$,

(3) ${\check\lambda}_{n-q_1k,-q_2k}=\gamma_m$ if $(n-q_1k,-q_2k)=(mq_1,mq_2)$ for some integer $m$.

The first case is impossible since $q_2k\neq0$. Let us demonstrate that the third case is also impossible. Indeed, in such a case
$$
n-q_1k=mq_1,\quad -q_2k=mq_2.
$$
from this $n=0$ that contradicts to the assumption $n\neq0$. Thus, the second case remains when
$$
n-q_1k=mp_1,\quad -q_2k=mp_2.
$$
From this
$$
k=-\frac{p_2}{p_1q_2-p_2q_1}\,n,\quad m=\frac{q_2}{p_1q_2-p_2q_1}\,n.
$$
Therefore the integers $p_2n$ and $q_2n$ should be divisible by $p_1q_2-p_2q_1$. Thus,
\begin{equation}
\sum\limits_{k\in\Z\setminus\{0\}}\gamma_k{\check\lambda}_{n-q_1k,-q_2k}=
\left\{\begin{array}{l}\beta_r\gamma_s\ \mbox{if}\ q_2n=r(p_1q_2-p_2q_1)\\
\qquad \mbox{and}\ p_2n=-s(p_1q_2-p_2q_1)\ (r,s\in\Z);\\
0\ \mbox{otherwise}.\end{array}\right.
                         \label{3.23}
\end{equation}

Right-hand sides of formulas \eqref{3.22} and \eqref{3.23} coincide. Therefore the left-hand side of the equation \eqref{3.21} is equal to zero for any $n\neq0$. Since $\alpha_n$ is not equal to zero at least for one $n\neq0$, we arrive to the second important conclusion: $a_1=0$. The equation \eqref{3.21} itself becomes a tautology and we forget it.

Now, we set $(n_1,n_2)=(np_1,np_2)$ in the equation \eqref{3.18}, where $n\neq0$. The pair $(n_1,n_2)$ numbers a node of the lattice $\Gamma'$ belonging to the line $L_1$. Taking the equalities $a_1=0$ and $(c_1,c_2)=(1,0)$ into account, the equation looks as follows:
\begin{equation}
\begin{aligned}
&np_2\sum\limits_{k\in\Z\setminus\{0\}}\alpha_k{\check\lambda}_{np_1-k,np_2}
+\varphi(np_1,np_2;p_1,p_2;q_1,q_2)\sum\limits_{k\in\Z\setminus\{0\}}\beta_k{\check\lambda}_{np_1-kp_1,np_2-kp_2}\\
&+\varphi(np_1,np_2;q_1,q_2;p_1,p_2)\sum\limits_{k\in\Z\setminus\{0\}}\gamma_k{\check\lambda}_{np_1-kq_1,np_2-kq_2}
=na_2p_2\beta_n.
\end{aligned}
                         \label{3.24}
\end{equation}
By \eqref{3.20},
$$
\varphi(np_1,np_2;p_1,p_2;q_1,q_2)=0,\quad \varphi(np_1,np_2;q_1,q_2;p_1,p_2)=-p_2n
$$
and the equation \eqref{3.24} is simplified to the following one:
\begin{equation}
\sum\limits_{k\in\Z\setminus\{0\}}\alpha_k{\check\lambda}_{np_1-k,np_2}
-\sum\limits_{k\in\Z\setminus\{0\}}\gamma_k{\check\lambda}_{np_1-kq_1,np_2-kq_2}
=a_2\beta_n.
                         \label{3.25}
\end{equation}

We first analyze the first sum on the left-hand side of \eqref{3.25}. Its summand $\alpha_k{\check\lambda}_{np_1-k,np_2}$ can be nonzero in three cases only:

(1) ${\check\lambda}_{np_1-k,np_2}=\alpha_m$ if $(np_1-k,np_2)=(m,0)$ for some integer $m$,

(2) ${\check\lambda}_{np_1-k,np_2}=\beta_m$ if $(np_1-k,np_2)=(mp_1,mp_2)$ for some integer $m$,

(3) ${\check\lambda}_{np_1-k,np_2}=\gamma_m$ if $(np_1-k,np_2)=(mq_1,mq_2)$ for some integer $m$.

The first case is impossible since $np_2\neq0$. Let us demonstrate that the second case is also impossible. Indeed, in such a case
$$
np_1-k=mp_1,\quad np_2=mp_2.
$$
From this $k=0$ that contradicts to the condition $k\in\Z\setminus\{0\}$. Thus, the third case remains when
$$
np_1-k=mq_1,\quad np_2=mq_2.
$$
From this
$$
k=\frac{p_1q_2-p_2q_1}{q_2}\,n,\quad m=\frac{p_2}{q_2}\,n.
$$
Hence the integer $p_2n$ must be divisible by $q_2$. If $p_2n=rq_2\ (r\in\Z)$, then $k=p_1n-q_1r$ and $(np_1-k,np_2)=(rq_1,rq_2)$. Thus,
\begin{equation}
\sum\limits_{k\in\Z\setminus\{0\}}\alpha_k{\check\lambda}_{np_1-k,np_2}=
\left\{\begin{array}{l}\alpha_{p_1n-q_1r}\gamma_r\ \mbox{if}\ p_2n=rq_2\ (r\in\Z\setminus\{0\}),\\
0\ \mbox{otherwise}.\end{array}\right.
                         \label{3.26}
\end{equation}

Next, we analyze the second sum on the left-hand side of \eqref{3.25}. Its summand\\ $\gamma_k{\check\lambda}_{np_1-kq_1,np_2-kq_2}$ can be nonzero in three cases only:

(1) ${\check\lambda}_{np_1-kq_1,np_2-kq_2}=\alpha_m$ if $(np_1-kq_1,np_2-kq_2)=(m,0)$ for some integer $m$,

(2) ${\check\lambda}_{np_1-kq_1,np_2-kq_2}=\beta_m$ if $(np_1-kq_1,np_2-kq_2)=(mp_1,mp_2)$ for some integer $m$,

(3) ${\check\lambda}_{np_1-kq_1,np_2-kq_2}=\gamma_m$ if $(np_1-kq_1,np_2-kq_2)=(mq_1,mq_2)$ for some integer $m$.

The second case is impossible. Indeed, in such a case the pair $(k,m)$ solves the system
$$
q_1k+p_1m=p_1n,\quad q_2k+p_2m=p_2n
$$
with nonzero determinant. The solution to the system is $k=0,m=n$. But $k\neq0$ in \eqref{3.25}.

The third case is also impossible. Indeed, in such a case the pair $(n,k+m)$ solves the homogeneous system
$$
p_1n-q_1(k+m)=0,\quad p_2n-q_2(k+m)=0
$$
with nonzero determinant. Hence $n=0$. But $n\neq0$ in \eqref{3.25}.

Thus, the first case remains. In this case we have the system
$$
np_1-kq_1=m,\quad np_2-kq_2=0.
$$
From this
$$
k=\frac{p_2}{q_2}n,\quad m=\frac{p_1q_2-p_2q_1}{q_2}n.
$$
Therefore $p_2n$ must be divisible by $q_2$. If $p_2n=rq_2\ (r\in\Z\setminus\{0\})$, then $k=r$ and $m=p_1n-q_1r$. Thus,
$$
\gamma_k{\check\lambda}_{np_1-kq_1,np_2-kq_2}=
\left\{\begin{array}{l}\gamma_r\alpha_{p_1n-q_1r},\ \mbox{if}\ p_2n=rq_2\ (r\in\Z\setminus\{0\}),\\
0\ \mbox{otherwise}.\end{array}\right.
$$
Hence
\begin{equation}
\sum\limits_{k\in\Z\setminus\{0\}}\gamma_k{\check\lambda}_{np_1-kq_1,np_2-kq_2}=
\left\{\begin{array}{l}\gamma_r\alpha_{p_1n-q_1r},\ \mbox{if}\ p_2n=rq_2\ (r\in\Z\setminus\{0\}),\\
0\ \mbox{otherwise}.\end{array}\right.
                         \label{3.27}
\end{equation}

Right-hand sides of formulas \eqref{3.26} and \eqref{3.27} coincide. Therefore the left-hand side of the equation \eqref{3.25} is equal to zero for any $n\neq0$. Since $\beta_n$ is not equal to zero at least for one $n\neq0$, we arrive to the conclusion: $a_2=0$. The equation \eqref{3.25} itself becomes a tautology and we forget it.

We have thus proved that
\begin{equation}
(a_1,a_2)=(0,0),\quad (c_1,c_2)=(1,0).
                         \label{3.28}
\end{equation}

In the same way we make sure with the help of \eqref{3.28} that the equation \eqref{3.18} is a tautology in the case when the pair $(n_1,n_2)$ numbers a node of the lattice $\Gamma'$ belonging to the line $L_2$, i.e., for $(n_1,n_2)=(nq_1,nq_2)\ (n\neq0)$.

\bigskip

Finally, we consider the equation \eqref{3.18} for a pair $(n_1,n_2)$ corresponding to a node of the lattice $\Gamma'$ which does not belong to
$L_0\cup L_1\cup L_2$, i.e., when
\begin{equation}
(n_1,n_2)\neq(0,0),\ (n_1,n_2)\neq(n,0),\ (n_1,n_2)\neq(np_1,np_2),\ (n_1,n_2)\neq(nq_1,nq_2).
                         \label{3.29}
\end{equation}
The right-hand side of \eqref{3.18} is equal to zero and the equation looks as follows:
\begin{equation}
\begin{aligned}
&n_2(p_1q_2-p_2q_1)\sum\limits_{k\in\Z\setminus\{0\}}\alpha_k{\check\lambda}_{n_1-k,n_2}
+q_2(p_2n_1-p_1n_2)\sum\limits_{k\in\Z\setminus\{0\}}\beta_k{\check\lambda}_{n_1-p_1k,n_2-p_2k}\\
&-p_2(q_2n_1-q_1n_2)\sum\limits_{k\in\Z\setminus\{0\}}\gamma_k{\check\lambda}_{n_1-q_1k,n_2-q_2k}=0.
\end{aligned}
                         \label{3.30}
\end{equation}
We have used \eqref{3.20} and $(c_1,c_2)=(1,0)$.

We first analyze the first sum on the left-hand side of \eqref{3.30}. By \eqref{3.29}, $n_2\neq0$. Therefore the summand of the sum $\alpha_k{\check\lambda}_{n_1-k,n_2}$ can be nonzero in two cases only:

(1) ${\check\lambda}_{n_1-k,n_2}=\beta_m$ if $(n_1-k,n_2)=(mp_1,mp_2)$ for some integer $m$,

(2) ${\check\lambda}_{n_1-k,n_2}=\gamma_m$ if $(n_1-k,n_2)=(mq_1,mq_2)$ for some integer $m$.

In the first case we have the system
$$
n_1-k=mp_1,\quad n_2=mp_2.
$$
Solving the system, we obtain $m=\frac{n_2}{p_2},\ n_1-k=n_1-\frac{p_1}{p_2}n_2$. Hence $n_2$ must be divisible by $p_2$. If $n_2=rp_2\ (r\in\Z\setminus\{0\})$, then
$m=r,\ n_1-k=rp_1$. Thus,
$$
\alpha_k{\check\lambda}_{n_1-k,n_2}=
\alpha_{n_1-rp_1}\beta_r,\ \mbox{if}\ n_2=rp_2\ \mbox{and}\ k=n_1-rp_1\ (r\in\Z\setminus\{0\}).
$$
The second case differs of the first one by the transposition $p\leftrightarrow q$ only. Therefore
$$
\alpha_k{\check\lambda}_{n_1-k,n_2}=
\alpha_{n_1-sq_1}\gamma_s,\ \mbox{if}\ n_2=sq_2\ \mbox{and}\ k=n_1-sq_1\ (s\in\Z\setminus\{0\}).
$$
It should be especially noted that the first and second cases take place simultaneously if $rp_2=sq_2$ for some $r,s\in\Z\setminus\{0\}$. But they occur for different values of $k$. Indeed, otherwise the pair $(r,s)$ would solve the linear homogeneous system
$$
rp_2=sq_2,\quad rp_1=sq_1
$$
with nonzero determinant. Hence $(r,s)=(0,0)$ that contradicts to the condition $n_2\neq0$.

From what was said in the previous paragraph it follows that
\begin{equation}
\sum\limits_{k\in\Z\setminus\{0\}}\alpha_k{\check\lambda}_{n_1-k,n_2}=
\left\{\begin{array}{l}\alpha_{n_1-rp_1}\beta_r+\alpha_{n_1-sq_1}\gamma_s\ \mbox{if}\ n_2=rp_2=sq_2\ (r,s\in\Z\setminus\{0\});\\
\alpha_{n_1-rp_1}\beta_r\ \mbox{if}\ n_2=rp_2\ (r\in\Z\setminus\{0\})\ \mbox{but}\ n_2\ \mbox{is not divisible by}\ q_2;\\
\alpha_{n_1-rq_1}\gamma_r\ \mbox{if}\ n_2=rq_2\ (r\in\Z\setminus\{0\})\ \mbox{but}\ n_2\ \mbox{is not divisible by}\ p_2;\\
0,\ \mbox{if}\ n_2\ \mbox{is not divisible either by}\ p_2\ \mbox{or by}\ q_2.\end{array}\right.
                         \label{3.33}
\end{equation}

Next, we analyze the second sum on the left-hand side of \eqref{3.30}. Its summand\\ $\beta_k{\check\lambda}_{n_1-p_1k,n_2-p_2k}$ can be nonzero in three cases only:

(1) ${\check\lambda}_{n_1-p_1k,n_2-p_2k}=\alpha_m$ if $(n_1-p_1k,n_2-p_2k)=(m,0)$ for some integer $m$,

(2) ${\check\lambda}_{n_1-p_1k,n_2-p_2k}=\beta_m$ if $(n_1-p_1k,n_2-p_2k)=(mp_1,mp_2)$ for some integer $m$,

(3) ${\check\lambda}_{n_1-p_1k,n_2-p_2k}=\gamma_m$ if $(n_1-p_1k,n_2-p_2k)=(mq_1,mq_2)$ for some integer $m$.

The second case is impossible. Indeed, in such a case
$$
n_1-p_1k=mp_1,\quad n_2-p_2k=mp_2.
$$
From this $(n_1,n_2)=\big((m+k)p_1,(m+k)p_2\big)$ that is prohibited by conditions \eqref{3.29}.

In the first case we have the system
$$
n_1-p_1k=m,\quad n_2-p_2k=0.
$$
From this $k=\frac{n_2}{p_2}$. Hence $n_2$ must be divisible by $p_2$. If $n_2=rp_2\ (r\in\Z)$, then
$k=r,\ m=n_1-rp_1$. Thus,
$$
\beta_k{\check\lambda}_{n_1-p_1k,n_2-p_2k}=
\beta_{r}\alpha_{n_1-p_1r}\ \mbox{if}\ n_2=rp_2\ \mbox{and}\ k=r\ (r\in\Z\setminus\{0\}).
$$
In the third case we have the system
$$
n_1-p_1k=mq_1,\quad n_2-p_2k=mq_2.
$$
From this
$$
k=\frac{q_2n_1-q_1n_2}{p_1q_2-p_2q_1},\quad m=\frac{-p_2n_1+p_1n_2}{p_1q_2-p_2q_1}.
$$
(Recall that $p_1q_2-p_2q_1\neq0$, see \eqref{3.8}.)
The integers $q_2n_1-q_1n_2$ and $-p_2n_1+p_1n_2$ must be divisible by $p_1q_2-p_2q_1$, i.e.,
$$
q_2n_1-q_1n_2=s(p_1q_2-p_2q_1),\quad -p_2n_1+p_1n_2=t(p_1q_2-p_2q_1)\quad(s,t\in\Z\setminus\{0\}).
$$
This is equivalent to the system
\begin{equation}
n_1=sp_1+tq_1,\quad n_2=sp_2+tq_2\quad(s,t\in\Z\setminus\{0\}).
                         \label{3.35}
\end{equation}
Hence
\begin{equation}
\beta_k{\check\lambda}_{n_1-p_1k,n_2-p_2k}=\left\{\begin{array}{l}
\beta_{s}\gamma_t\ \mbox{if}\ n_1=sp_1+tq_1, n_2=sp_2+tq_2\ (s,t\in\Z\setminus\{0\})\ \mbox{and}\ k=s;\\
0\ \mbox{otherwise}.\end{array}\right.
                         \label{3.36}
\end{equation}

First and third cases have a nonempty intersection that is characterized by the relations
\begin{equation}
n_1=sp_1+tq_1,\quad n_2=rp_2=sp_2+tq_2\quad(r,s,t\in\Z\setminus\{0\}).
                         \label{3.37}
\end{equation}
From this $(r-s)p_2=tq_2$, i.e., $tq_2$ is divisible by $p_2$. Conversely, if the relations \eqref{3.36} hold and $tq_2$ is divisible by $p_2$, then the first case takes place.

The situation when the first case takes place but the third case does not take is characterized by the following:
\begin{equation}
\begin{aligned}
&n_2=rp_2\ (r\in\Z\setminus\{0\}),\ \mbox{but the representation}\\
&n_1=sp_1+tq_1,\ n_2=sp_2+tq_2\ (s,t\in\Z\setminus\{0\})\ \mbox{is impossible}.
\end{aligned}
                         \label{3.38}
\end{equation}

The situation when the third case takes place but the first case does not take is characterized by the following:
\begin{equation}
n_1=sp_1+tq_1,\ n_2=sp_2+tq_2\ (s,t\in\Z\setminus\{0\})\ \mbox{and}\ n_2\ \mbox{is not divisible by}\ p_2.
                         \label{3.39}
\end{equation}

Now, we compute the sum $\sum\limits_{k\in\Z\setminus\{0\}}\beta_k{\check\lambda}_{n_1-p_1k,n_2-p_2k}$. It looks different in four cases:

(a) If relations \eqref{3.37} hold, then, by \eqref{3.35} and \eqref{3.36}, the sum has only two nonzero summands corresponding to $k=r$ and $k=s$. These summands are different, i.e., $r\neq s$ as is seen from \eqref{3.37}. Thus,
\begin{equation}
\begin{aligned}
\sum\limits_{k\in\Z\setminus\{0\}}&\beta_k{\check\lambda}_{n_1-p_1k,n_2-p_2k}=
\alpha_{n_1-p_1r}\beta_{r}+\beta_s\gamma_t,\\
&\mbox{if}\ n_1=sp_1+tq_1,\quad n_2=rp_2=sp_2+tq_2\quad(r,s,t\in\Z\setminus\{0\}).
\end{aligned}
                         \label{3.40}
\end{equation}

(b) If relations \eqref{3.38} hold, then the sum $\sum\limits_{k\in\Z\setminus\{0\}}\beta_k{\check\lambda}_{n_1-p_1k,n_2-p_2k}$ contains only one nonzero summand corresponding to $k=r$. Thus,
\begin{equation}
\begin{aligned}
\sum\limits_{k\in\Z\setminus\{0\}}&\beta_k{\check\lambda}_{n_1-p_1k,n_2-p_2k}=
\alpha_{n_1-p_1r}\beta_{r},\\
&\mbox{if}\ n_2=rp_2\ (r\in\Z\setminus\{0\})\ \mbox{but the representation}\\
&n_1=sp_1+tq_1,\ n_2=sp_2+tq_2\ (s,t\in\Z\setminus\{0\})\ \mbox{is impossible}.
\end{aligned}
                         \label{3.41}
\end{equation}

(c) If relations \eqref{3.39} hold, then the sum $\sum\limits_{k\in\Z\setminus\{0\}}\beta_k{\check\lambda}_{n_1-p_1k,n_2-p_2k}$ has only one nonzero summand corresponding to $k=s$. Thus,
\begin{equation}
\begin{aligned}
&\sum\limits_{k\in\Z\setminus\{0\}}\beta_k{\check\lambda}_{n_1-p_1k,n_2-p_2k}=
\beta_s\gamma_t,\\
&\mbox{if}\ n_1=sp_1+tq_1,\ n_2=sp_2+tq_2\ (s,t\in\Z\setminus\{0\})\ \mbox{and}\ n_2\ \mbox{is not divisible by}\ p_2.
\end{aligned}
                         \label{3.42}
\end{equation}

(d) Finally, the situation can occur when $n_2$ is not divisible by $p_2$ and the representation
$n_1=sp_1+tq_1,\ n_2=sp_2+tq_2\ (s,t\in\Z\setminus\{0\})$ is impossible. In such a case none of relations \eqref{3.37}--\eqref{3.39} are fulfill and $\sum\limits_{k\in\Z\setminus\{0\}}\beta_k{\check\lambda}_{n_1-p_1k,n_2-p_2k}=0$.

We unite formulas \eqref{3.40}--\eqref{3.42}
\begin{equation}
\sum\limits_{k\in\Z\setminus\{0\}}\beta_k{\check\lambda}_{n_1-p_1k,n_2-p_2k}=
\left\{\begin{array}{l}
\alpha_{n_1-p_1r}\beta_{r}+\beta_s\gamma_t\ \mbox{if}\\
\quad n_1=sp_1+tq_1,\ n_2=rp_2=sp_2+tq_2\ (r,s,t\in\Z\setminus\{0\});\\
\alpha_{n_1-p_1r}\beta_{r}\ \mbox{if}\ n_2=rp_2\ (r\in\Z\setminus\{0\})\
\mbox{but the representation}\\
\quad n_1=sp_1+tq_1,\ n_2=sp_2+tq_2\ (s,t\in\Z\setminus\{0\})\ \mbox{is impossible};\\
\beta_s\gamma_t\
\mbox{if}\ n_1=sp_1+tq_1,\ n_2=sp_2+tq_2\ (s,t\in\Z\setminus\{0\})\\
\qquad \mbox{and}\ n_2\ \mbox{is not divisible by}\ p_2;\\
0\ \mbox{if}\ n_2\ \mbox{is not divisible by}\ p_2\ \mbox{and the representation}\\
\quad n_1=sp_1+tq_1,\ n_2=sp_2+tq_2\ (s,t\in\Z\setminus\{0\})\ \mbox{is impossible}.
\end{array}\right.
                         \label{3.43}
\end{equation}

The third sum on the left-hand side of \eqref{3.30} is obtained from the second sum by the transposition $\beta\leftrightarrow\gamma,\,p\leftrightarrow q$. Therefore we just perform the transposition $\beta\leftrightarrow\gamma,\,p\leftrightarrow q,\,s\rightarrow t$ in \eqref{3.43} to obtain
\begin{equation}
\sum\limits_{k\in\Z\setminus\{0\}}\gamma_k{\check\lambda}_{n_1-q_1k,n_2-q_2k}=
\left\{\begin{array}{l}
\alpha_{n_1-q_1r}\gamma_{r}+\beta_s\gamma_t\ \mbox{if}\\
\quad n_1=sp_1+tq_1,\ n_2=rq_2=sp_2+tq_2\ (r,s,t\in\Z\setminus\{0\});\\
\alpha_{n_1-q_1r}\gamma_{r}\ \mbox{if}\ n_2=rq_2\ (r\in\Z\setminus\{0\})\
\mbox{but the representation}\\
\quad n_1=sp_1+tq_1,\ n_2=sp_2+tq_2\ (s,t\in\Z\setminus\{0\})\ \mbox{is impossible};\\
\beta_s\gamma_t\
\mbox{if}\ n_1=sp_1+tq_1,\ n_2=sp_2+tq_2\ (s,t\in\Z\setminus\{0\})\\
\quad \mbox{and}\ n_2\ \mbox{is not divisible by}\ q_2;\\
0\ \mbox{if}\ n_2\ \mbox{is not divisible by}\ q_2\ \mbox{and the representation}\\
\quad n_1=sp_1+tq_1,\ n_2=sp_2+tq_2\ (s,t\in\Z\setminus\{0\})\ \mbox{is impossible}.
\end{array}\right.
                         \label{3.44}
\end{equation}

\bigskip

Now, we have to substitute expressions \eqref{3.33} and \eqref{3.43}--\eqref{3.44} into the equation \eqref{3.30}. The substitution is not easy since each of the expressions consists of four lines. Formally speaking, we have to consider $4\times4\times4=64$ combinations. Fortunately, the most of these combinations is logically impossible since they contain conditions contradicting to each other. An easy but bulky analysis shows that there exist 7 logically possible versions of the equation \eqref{3.30}. For our purposes, it suffices to consider the following three versions.

Version 1. Assume that $n_2$ is not divisible either by $p_2$ or by $q_2$ and the representation $n_1=sp_1+tq_1,\ n_2=sp_2+tq_2\ (s,t\in\Z\setminus\{0\})$ is possible. In such a case formulas \eqref{3.33} and \eqref{3.43}--\eqref{3.44} give
$$
\sum\limits_{k\in\Z\setminus\{0\}}\alpha_k{\check\lambda}_{n_1-k,n_2}=0,\quad
\sum\limits_{k\in\Z\setminus\{0\}}\beta_k{\check\lambda}_{n_1-p_1k,n_2-p_2k}=\beta_s\gamma_t,\quad
\sum\limits_{k\in\Z\setminus\{0\}}\gamma_k{\check\lambda}_{n_1-q_1k,n_2-q_2k}=\beta_s\gamma_t.
$$
Substituting these expressions into the equation \eqref{3.30}, we obtain
$$
n_2(p_1q_2-p_2q_1)\beta_s\gamma_t=0.
$$
Since $n_2\neq0$ and $p_1q_2-p_2q_1\neq0$, we conclude $\beta_s\gamma_t=0$. Thus,
\begin{equation}
\left.
\begin{aligned}
\beta_s\gamma_t=0\ &\mbox{if}\ n_2\ \mbox{is not divisible either by}\ p_2\ \mbox{or by}\ q_2\\
&\mbox{and}\  n_1=sp_1+tq_1,\ n_2=sp_2+tq_2\ (s,t\in\Z\setminus\{0\}).
\end{aligned}
\right\}
                         \label{3.45}
\end{equation}

Version 2. Assume that
$$
n_2=rp_2=uq_2,\quad n_1-rp_2\neq0,\quad n_1-uq_2\neq0\quad(r,u\in\Z\setminus\{0\})
$$
and the representation $n_1=sp_1+tq_1,\ n_2=sp_2+tq_2\ (s,t\in\Z\setminus\{0\})$ is impossible. In such a case formulas \eqref{3.33}, \eqref{3.43}--\eqref{3.44} give
$$
\begin{aligned}
\sum\limits_{k\in\Z\setminus\{0\}}\alpha_k{\check\lambda}_{n_1-k,n_2}&=\alpha_{n_1-rp_1}\beta_r+\alpha_{n_1-uq_1}\gamma_u,\\
\sum\limits_{k\in\Z\setminus\{0\}}\beta_k{\check\lambda}_{n_1-p_1k,n_2-p_2k}&=\alpha_{n_1-p_1r}\beta_r,\quad
\sum\limits_{k\in\Z\setminus\{0\}}\gamma_k{\check\lambda}_{n_1-q_1k,n_2-q_2k}=\alpha_{n_1-q_1u}\gamma_u.
\end{aligned}
$$
Substituting these values into \eqref{3.30}, we obtain
\begin{equation}
\left.
\begin{aligned}
&p_2(q_2n_1-q_1n_2)\alpha_{n_1-rp_1}\beta_r-q_2(p_2n_1-p_1n_2)\alpha_{n_1-uq_1}\gamma_u=0\\
&\quad\mbox{if}\ n_2=rp_2=uq_2,\ n_1-rp_1\neq0, n_1-uq_1\neq0,\ (r,u\in\Z\setminus\{0\})\ \mbox{and the}\\
&\quad \mbox{representation}\ n_1=sp_1+tq_1,\ n_2=sp_2+tq_2\ (s,t\in\Z\setminus\{0\})\ \mbox{is impossible}.
\end{aligned}
\right\}
                         \label{3.46}
\end{equation}

Version 3. Assume that
$$
n_2=rp_2=uq_2,\quad n_1-rp_2\neq0,\quad n_1-uq_2\neq0\quad(r,u\in\Z\setminus\{0\})
$$
and the representation $n_1=sp_1+tq_1,\ n_2=sp_2+tq_2\ (s,t\in\Z\setminus\{0\})$ is possible. In such a case formulas \eqref{3.33}, \eqref{3.43}--\eqref{3.44} give
$$
\begin{aligned}
\sum\limits_{k\in\Z\setminus\{0\}}\alpha_k{\check\lambda}_{n_1-k,n_2}&=\alpha_{n_1-rp_1}\beta_r+\alpha_{n_1-uq_1}\gamma_u,\\
\sum\limits_{k\in\Z\setminus\{0\}}\beta_k{\check\lambda}_{n_1-p_1k,n_2-p_2k}&=\alpha_{n_1-rp_1}\beta_r+\beta_s\gamma_t,\quad
\sum\limits_{k\in\Z\setminus\{0\}}\gamma_k{\check\lambda}_{n_1-q_1k,n_2-q_2k}=\alpha_{n_1-uq_1}\gamma_u+\beta_s\gamma_t.
\end{aligned}
$$
Substituting these values into the equation \eqref{3.30}, we obtain
\begin{equation}
\left.
\begin{aligned}
p_2(q_2n_1-q_1n_2)&\alpha_{n_1-rp_1}\beta_r-q_2(p_2n_1-p_1n_2)\alpha_{n_1-uq_1}\gamma_u-n_2(p_1q_2-p_2q_1)\beta_{s}\gamma_t=0\\
&\mbox{if}\ n_2=rp_2=uq_2,\ n_1-rp_1\neq0, n_1-uq_1\neq0\\
&\mbox{and}\   n_1=sp_1+tq_1,\ n_2=sp_2+tq_2\quad(r,s,t,u\in\Z\setminus\{0\});
\end{aligned}
\right\}
                         \label{3.47}
\end{equation}

We emphasize that \eqref{3.45}--\eqref{3.47} are three versions of the same equation \eqref{3.30} and these versions do not exhaust the latter equation (4 other versions exist but we do not use them).

Now, we are going to transform equations \eqref{3.45}--\eqref{3.47} by simplifying the involved conditions.
Recall that $(p_2,q_2)$ are positive integers. Let $d>0$ be their greatest common divisor. Then $(p_2,q_2)=(dp'_2,dq'_2)$, where $(p'_2,q'_2)$ is the pair of relatively prime positive integers. Now $p_1q_2-p_2q_1=d(p_1q'_2-p'_2q_1)$. Let us introduce the integer $\delta=p_1q'_2-p'_2q_1$, it will participate in many further formulas.

We start with the equation \eqref{3.45}.
Observe that $n_1$ and $n_2$ do not participate in the equation itself but are involved to conditions written in the second line of \eqref{3.45}. For arbitrary $s,t\in \Z\setminus\{0\}$, we can define $n_1$ and $n_2$ by the equalities
$$
n_1=sp_1+tq_1,\quad n_2=sp_2+tq_2.
$$
It remains to ensure compliance with the condition: $n_2$ is not divisible by either $p_2$ or by $q_2$. Writing the second equality in the form 
$ n_2=dsp'_2+dtq'_2$, we see that the condition is equivalent to the following one:
$$
s\neq0\ \mbox{is not divisible by}\ q'_2\ \mbox{and}\ t\neq0\ \mbox{is not divisible by}\ p'_2.
$$
Thus, the equation \eqref{3.45} takes the form
$$
\beta_s\gamma_t=0\
\mbox{if}\ s\neq0\ \mbox{is not divisible by}\ q'_2\ \mbox{and}\ t\neq0\ \mbox{is not divisible by}\ p'_2
$$
or, after changing notations,
\begin{equation}
\beta_{n_1}\gamma_{n_2}=0\
\mbox{if}\ n_1\neq0\ \mbox{is not divisible by}\ q'_2\ \mbox{and}\ n_2\neq0\ \mbox{is not divisible by}\ p'_2.
                         \label{3.48}
\end{equation}

Next, we transform the equation \eqref{3.46} that takes the form
\begin{equation}
\left.
\begin{aligned}
&p'_2(dq'_2n_1-q_1n_2)\alpha_{n_1-rp_1}\beta_r-q'_2(dp'_2n_1-p_1n_2)\alpha_{n_1-uq_1}\gamma_u=0\\
&\quad\mbox{if}\ n_2=rp_2=uq_2,\ n_1-rp_1\neq0, n_1-uq_1\neq0,\ (r,u\in\Z\setminus\{0\})\ \mbox{and the}\\
&\quad\mbox{representation}
\ n_1=sp_1+tq_1,\ n_2=sp_2+tq_2\ (s,t\in\Z\setminus\{0\})\ \mbox{is impossible}.
\end{aligned}
\right\}
                         \label{3.49}
\end{equation}
Since $p'_2$ and $q'_2$ are relatively prime integers, the equalities $n_2=rdp'_2=udq'_2$ imply the existence of the integer $k$ such that
\begin{equation}
n_2=dp'_2q'_2k,\quad r=q'_2k,\quad u=p'_2k,
                         \label{3.50}
\end{equation}
and the equation \eqref{3.49} takes the form
\begin{equation}
(n_1-p'_2q_1k)\alpha_{n_1-p_1q'_2k}\beta_{q'_2k}-(n_1-p_1q'_2k)\alpha_{n_1-p'_2q_1k}\gamma_{p'_2k}=0.
                         \label{3.51}
\end{equation}
The integer $n=n_1-p_1q'_2k$ is as arbitrary as $n_1$. Expressing $n_1-p'_2q_1k=n+\delta k$ and substituting these expressions into \eqref{3.51}, we write the result in the form
$$
\frac{\alpha_{n}}{n}\,\frac{\beta_{q'_2k}}{k}-\frac{\alpha_{n+\delta k}}{n+\delta k}\,\frac{\gamma_{p'_2k}}{k}=0\quad(n\neq0,k\neq0,n+\delta k\neq0).
$$

We also have to discuss conditions written in last two lines of \eqref{3.49}. Taking the equality $n_2=dp'_2q'_2k$ into account, the system
$n_1=sp_1+tq_1,\ n_2=sp_2+tq_2$ can be written as follows:
$$
p_1s+q_1t=n_1,\quad p'_2s+q'_2t=p'_2q'_2k.
$$
From this
$$
s=\frac{q'_2(n_1-p'_2q_1k)}{\delta},\quad t=\frac{p'_2(-n_1+p_1q'_2k)}{\delta}.
$$
Using the equalities $n_1-p'_2q_1k=n+\delta k,\ n_1-p_1q'_2k=n$, we write the latter formula in the form
$$
s=\frac{q'_2n}{\delta}+q'_2k,\quad t=-\frac{p'_2n}{\delta}.
$$
By conditions written in last two lines of \eqref{3.49}, at least one of reals $s,t$ is not an integer. Hence these conditions are equivalent to the following one:
\begin{equation}
(p'_2n\ \mbox{is not divisible by}\ \delta)\ \mbox{or}\ (q'_2n\ \mbox{is not divisible by}\ \delta).
                         \label{3.53}
\end{equation}
Let us demonstrate that this is equivalent to the condition: $n$ is not divisible by $\delta$.
Indeed, the negation of \eqref{3.53} means that
$$
p'_2n=\delta a,\quad q'_2n=\delta b
$$
with some integers $a$ and $b$. Taking a linear combination of these equalities, we have
$$
(lp'_2+mq'_2)n=\delta (la+mb).
$$
Since $p'_2$ and $q'_2$ are relatively prime integers, there exist integers $l$ and $m$ such that $lp'_2+mq'_2=1$. Now, the last formula gives $n=\delta (la+mb)$, i.e., $n$ is divisible by $\delta$.

Thus, the complete form of the equation \eqref{3.49} looks as follows:
$$
\frac{\alpha_{n}}{n}\,\frac{\beta_{q'_2k}}{k}-\frac{\alpha_{n+\delta k}}{n+\delta k}\,\frac{\gamma_{p'_2k}}{k}=0\
\mbox{if}\ n\neq0, k\neq0\ \mbox{and}\ n\ \mbox{is not divisible by}\ \delta.
$$
Changing notation of indices, we write this in the final form
\begin{equation}
\frac{\alpha_{n_1}}{n_1}\,\frac{\beta_{q'_2n_2}}{n_2}-\frac{\alpha_{n_1+\delta n_2}}{n_1+\delta n_2}\,\frac{\gamma_{p'_2n_2}}{n_2}=0\
\mbox{if}\ n_1\neq0, n_2\neq0\ \mbox{and}\ n_1\ \mbox{is not divisible by}\ \delta.
                         \label{3.55}
\end{equation}

Finally, we transform the equation \eqref{3.47} that takes the form
\begin{equation}
\left.
\begin{aligned}
&p'_2(dq'_2n_1-q_1n_2)\alpha_{n_1-rp_1}\beta_r-q'_2(dp'_2n_1-p_1n_2)\alpha_{n_1-uq_1}\gamma_u-n_2\delta\,\beta_{s}\gamma_t=0\\
&\quad\mbox{if}\ n_2=rdp'_2=udq'_2,\ n_1-rp_1\neq0, n_1-uq_1\neq0,\\
&\quad n_1=sp_1+tq_1,\ n_2=d(sp'_2+tq'_2)\quad(r,s,t,u\in\Z\setminus\{0\});
\end{aligned}
\right\}
                         \label{3.56}
\end{equation}
Since $p'_2$ and $q'_2$ are relatively prime integers, the equalities $n_2=rdp'_2=udq'_2$ imply the validity of \eqref{3.50} with some integer $k$.
Now, the system
$n_1=sp_1+tq_1,\ n_2=d(sp'_2+tq'_2)$ takes the form
\begin{equation}
p_1s+q_1t=n_1,\quad p'_2s+q'_2t=p'_2q'_2k.
                         \label{3.57}
\end{equation}
By the condition, this system has the integer solution $(s,t)$. Hence $s$ is divisible by $q'_2$ and $t$ is divisible by $p'_2$: $s=q'_2s',\ t=p'_2t'$. The system takes the form
$$
p_1q'_2s'+p'_2q_1t'=n_1,\quad s'+t'=k.
$$
From this
$$
s'=\frac{n_1-p'_2q_1k}{\delta},\quad t'=\frac{-n_1+p_1q'_2k}{\delta}.
$$
Hence the integers $n_1-p'_2q_1k$ and $n_1-p_1q'_2k$ are divisible by $\delta$, i.e.,
\begin{equation}
n_1=\delta s'+p'_2q_1k=-\delta t'+p_1q'_2k.
                         \label{3.58}
\end{equation}

Equations \eqref{3.50} and \eqref{3.58} imply
\begin{equation}
n_1-rp_1=-\delta t',\quad n_1-uq_1=\delta s'.
                         \label{3.59}
\end{equation}
Substituting these values and $n_2=dp'_2q'_2k$ into \eqref{3.56}, we obtain
$$
(n_1-p'_2q_1k)\alpha_{-\delta t'}\beta_{q'_2k}-(n_1-p_1q'_2k)\alpha_{\delta s'}\gamma_{p'_2k}-\delta k\,\beta_{s}\gamma_t=0.
$$
In virtue of \eqref{3.50} formulas \eqref{3.59} can be written as $n_1-p'_2q_1k =\delta s'$, $n_1-p_1q'_2k=-\delta t'$.
This allows us to simplify the latter equation to the following one:
$$
s'\,\alpha_{-\delta t'}\beta_{q'_2k}+t'\,\alpha_{\delta s'}\gamma_{p'_2k}-k\,\beta_{s}\gamma_t=0.
$$
Finally, inserting the values $s=q'_2 s',\ t=p'_2 t',\ k=s'+t'$, we obtain
$$
s'\,\alpha_{-\delta t'}\beta_{q'_2(s'+t')}+t'\,\alpha_{\delta s'}\gamma_{p'_2(s'+t')}-(s'+t')\,\beta_{q'_2 s'}\gamma_{p'_2 t'}=0.
$$
Changing indices according to the equalities $n_1=-t',\ n_2=s'+t'$, we arrive to the equation
$$
(n_1+n_2)\alpha_{\delta n_1}\beta_{q'_2n_2}-n_1\alpha_{\delta(n_1+n_2)}\gamma_{p'_2n_2}-n_2\beta_{q'_2(n_1+n_2)}\gamma_{-p'_2n_1}=0\quad
(n_1\neq0,n_2\neq0,n_1+n_2\neq0)
$$
that can be written in the final form
\begin{equation}
\frac{\alpha_{\delta n_1}}{n_1}\,\frac{\beta_{q'_2n_2}}{n_2}-
\frac{\alpha_{\delta(n_1+n_2)}}{n_1+n_2}\frac{\gamma_{p'_2n_2}}{n_2}
-\frac{\beta_{q'_2(n_1+n_2)}}{n_1+n_2}\,\frac{\gamma_{-p'_2n_1}}{n_1}=0\quad
(n_1\neq0,n_2\neq0,n_1+n_2\neq0),
                         \label{3.60}
\end{equation}

\bigskip

We can now finish the proof of Theorem \ref{Th1.2}.
Introducing the variables
\begin{equation}
x_n=\frac{\alpha_{\delta n}}{n},\quad y_n=\frac{\beta_{q'_2 n}}{n},\quad z_n=\frac{\gamma_{p'_2n}}{n}\quad(n\neq0),
                         \label{3.61}
\end{equation}
we write the equation \eqref{3.60} in the form
$$
x_{n_1}y_{n_2}-x_{n_1+n_2}z_{n_2}+y_{n_1+n_2}z_{-n_1}=0\quad
(n_1\neq0,n_2\neq0,n_1+n_2\neq0)
$$
that coincides with \eqref{1.13}. Hypotheses \eqref{1.14} and \eqref{1.15} of Theorem \ref{Th1.3} are satisfied by \eqref{3.11} and \eqref{3.12}. Applying Theorem \ref{Th1.3}, we can state that one of the following three cases takes place: (1) $y_n=z_n=0$ for all $n\neq0$, (2) $x_n=z_n=0$ for all $n\neq0$, (3) $x_n=y_n=0$ for all $n\neq0$.

Recall that we are proving Theorem \ref{Th1.2} by contradiction. At the beginning of the current section, we assumed that each of lines \eqref{3.2} contained at least one point belonging to the spectrum of the function $\lambda$ and other than the origin. This means that each of the sequences
$\{\alpha_n\}_{n\in\Z\setminus\{0\}}, \{\beta_n\}_{n\in\Z\setminus\{0\}}, \{\gamma_n\}_{n\in\Z\setminus\{0\}}$ contains at least one nonzero term. We will get a contradiction in each of the cases presented by Theorem \ref{Th1.3}.

(1) $y_n=z_n=0$ for all $n\neq0$. By \eqref{3.61}, this means that
\begin{equation}
\beta_{q'_2n}=\gamma_{p'_2n}=0\quad (n\neq0).
                         \label{3.63}
\end{equation}
We choose $m,l\in\Z\setminus\{0\}$ such that $\beta_m\neq0$ and $\gamma_l\neq0$. By \eqref{3.63}, $m$ is not divisible by $q'_2$ and $l$ is not divisible by $p'_2$. Then  $\beta_m\gamma_l\neq0$. But by the equation \eqref{3.48}, $\beta_m\gamma_l=0$. We have got a contradiction.

(2) $x_n=z_n=0$ for all $n\neq0$. By \eqref{3.61}, this means that
\begin{equation}
\alpha_{\delta n}=\gamma_{p'_2n}=0\quad (n\neq0).
                         \label{3.64}
\end{equation}
The equation \eqref{3.55} is simplified to the following one:
$$
\alpha_{n_1}\beta_{q'_2n_2}=0\quad
\mbox{if}\ n_1\neq0, n_2\neq0\ \mbox{and}\ n_1\ \mbox{is not divisible by}\ \delta.
$$
If we assumed that $\beta_{q'_2n_2}\neq0$ at least for one $n_2\neq0$, then the last formula would imply that
$$
\alpha_{n}=0\quad\mbox{for all}\ n\neq0,\ \mbox{that are not divisible by}\ \delta.
$$
Together with \eqref{3.64}, this gives the contradiction: $\alpha_n=0$ for all $n\neq0$. Thus, only one possibility remains: $\beta_{q'_2n}=0$ for all $n\neq0$. Together with \eqref{3.64}, this gives $y_n=z_n=0\ (n\neq0)$ and we return to the case (1) considered above.

The third case $x_n=y_n=0\ (n\neq0)$ is considered in the same way.

\section{Proof of Theorem \ref{Th1.3}}

We investigate the system \eqref{1.13} together with the parity condition \eqref{1.14}. The latter condition will be multiply used as a ``default condition'', i.e., with no reference.

We first list some obvious but important properties of the system \eqref{1.13}. First of all the system possesses the following homogeneity: if $(x_n,y_n,z_n)$ is a solution to the system, then $(ax_n,ay_n,az_n)$ is also a solution for any real $a$.
Second, the system is invariant under the changes
\begin{equation}
\begin{aligned}
&(n_1,n_2)\rightarrow(-n_1,-n_2);\\
&(x_n,y_n,z_n)\rightarrow (y_n,x_n,-z_{-n}),\quad(n_1,n_2)\rightarrow(n_2,n_1);\\
&(x_n,y_n,z_n)\rightarrow (z_n,y_n,x_{n}),\quad(n_1,n_2)\rightarrow(-n_1,n_1+n_2);\\
&(x_n,y_n,z_n)\rightarrow (-x_{-n},z_n,y_{n}),\quad(n_1,n_2)\rightarrow(-n_1-n_2,n_2).
\end{aligned}
                         \label{4.1}
\end{equation}

Obviously, the system has one-dimensional solutions: if $y_n=z_n=0$ for all $n\neq0$, then $x_n$ can be arbitrary; the same is true in two other cases: $x_n=z_n=0$ and $x_n=y_n=0$. We are going to prove that the system has no other solution. To this end the decay condition \eqref{1.15} will be used at a couple of crucial points.

Let us first prove that the system has no two-dimensional solution; more precisely: every two-dimensional solution is actually a one-dimensional one. Indeed, let $z_n=0$ for all $n\neq0$. The equation \eqref{1.13} is simplified to the following one:
$$
x_{n_1}y_{n_2}=0\quad(n_1\neq0,n_2\neq0,n_1+n_2\neq0).
$$
If we assumed that $x_{n_1^0}\neq0$ for some $n_1^0\neq0$, then the last formula would imply that $y_{n_2}=0$ for all $n_2\neq0,\ n_2\neq-n_1^0$. In particular,
$y_{n_1^0}=0$ and hence $y_{-n_1^0}=\overline{y_{n_1^0}}=0$. Thus, $y_n=z_n=0$ for all $n\neq0$, i.e., our solution is one-dimensional. Two other cases $y_n=0$ and $x_n=0$ are considered in the same way.

The left-hand side of the equation \eqref{1.13} is a quadratic form in the variables $(x,y,z)$ but it is linear in each of the variables. Such forms are called trilinear forms.

For every $n\ge2$, the system \eqref{1.13} contains three finite subsystems:
\begin{equation}
z_{n-m}x_{n}-z_{-m}y_{n}=x_{m}y_{n-m}\quad(1\le m\le n-1),
                         \label{4.2}
\end{equation}
\begin{equation}
x_{-m}y_{n}-x_{n-m}z_{n}=-y_{n-m}z_{m}\quad(1\le m\le n-1),
                         \label{4.3}
\end{equation}
\begin{equation}
y_{-m}x_{n}+y_{n-m}z_{-n}=x_{n-m}z_{-m}\quad(1\le m\le n-1).
                         \label{4.4}
\end{equation}
If $x_k,y_k,z_k$ are assumed to be known for $|k|<n$, then \eqref{4.2} can be considered as a system of linear equations with two unknowns $(x_n,y_n)$. The system is overdetermined for $n\ge3$ and the degree of overdetermination grows with $n$. The same is true for systems \eqref{4.3} and \eqref{4.4}.

We now derive an interesting corollary of equations \eqref{4.2}--\eqref{4.4}. To this end we write down equations \eqref{4.2} and \eqref{4.4} together, performing the change $m:=n-m$ in the latter equation
$$
\begin{array}{c}
z_{n-m}x_{n}-z_{-m}y_{n}=x_{m}y_{n-m},\\
y_{-n+m}x_{n}+y_{m}z_{-n}=x_{m}z_{-n+m}.
\end{array}
\left|
\begin{array}{c}z_{-n+m}\\-y_{n-m}\end{array}\right.
$$
We eliminate $x_{m}$ from this system by multiplying the first equation by $z_{-n+m}$, multiplying the second equation by $-y_{n-m}$ (as is designated after the vertical bar) and summing the results. In this way we obtain
$$
(|z_{n-m}|^2-|y_{n-m}|^2)x_{n}=z_{-n+m}z_{-m}y_{n}+y_{m}y_{n-m}z_{-n}.
$$
Let us replace $m$ with $n-m$ here
$$
(|z_{m}|^2-|y_{m}|^2)x_{n}=z_{-n+m}z_{-m}y_{n}+y_{m}y_{n-m}z_{-n}.
$$
Right-hand sides of two last equations coincide. Equating left-hand sides, we get
$$
(|z_{m}|^2-|y_{m}|^2)x_{n}=(|z_{n-m}|^2-|y_{n-m}|^2)x_{n}\quad(1\le m\le n-1).
$$

The same trick can be done with pairs \eqref{4.2}--\eqref{4.3} and \eqref{4.3}--\eqref{4.4}. We thus obtain
\begin{equation}
\left.
\begin{aligned}
(|z_{m}|^2-|y_{m}|^2)x_{n}&=(|z_{n-m}|^2-|y_{n-m}|^2)x_{n},\\
(|x_{m}|^2-|z_{m}|^2)y_{n}&=(|x_{n-m}|^2-|z_{n-m}|^2)y_{n},\\
(|x_{m}|^2-|y_{m}|^2)z_{n}&=(|x_{n-m}|^2-|y_{n-m}|^2)z_{n}.
\end{aligned}
\right\}\quad(1\le m\le n-1).
                         \label{4.5}
\end{equation}

In this section, we will show that Theorem \ref{Th1.3} follows from its partial case.

\begin{lemma} \label{L4.1}
Theorem \ref{Th1.3} is valid under the additional condition
$$
|x_1|+|y_1|+|z_1|>0.
$$
\end{lemma}

We start the proof of Theorem \ref{Th1.3}. Let a sequence $(x_n,y_n,z_n)_{n\in\Z\setminus\{0\}}$ satisfy hypotheses of the theorem. We can assume that the sequence is not identically equal to zero (otherwise there is nothing to prove). Let $ p$ be the minimal of positive integers $n$ such that $|x_n|+|y_n|+|z_n|>0$.
If $p=1$, we are under hypotheses of Lemma \ref{L4.1} that is assumed to be valid. Therefore we can assume that $p\ge2$.
The sequence
$({\tilde x}_n,{\tilde y}_n,{\tilde z}_n)=(x_{ p n},y_{ p n},z_{ p n})$ also satisfy hypotheses of Theorem \ref{Th1.3} and
$|{\tilde x}_1|+|{\tilde y}_1|+|{\tilde z}_1|>0$. By Lemma \ref{L4.1}, the assertion of Theorem \ref{Th1.3} is true for $({\tilde x}_n,{\tilde y}_n,{\tilde z}_n)$, i.e., one of the following three statements is valid: (1) ${\tilde y}_n={\tilde z}_n=0\ (n\neq0)$, (2) ${\tilde x}_n={\tilde z}_n=0\ (n\neq0)$, (3) ${\tilde x}_n={\tilde y}_n=0\ (n\neq0)$. In view of symmetries \eqref{4.1}, we can assume without lost of generality that the first statement is valid. Thus,
\begin{equation}
y_{ p n}=z_{ p n}=0\quad (n\in\Z\setminus\{0\}).
                         \label{4.6}
\end{equation}
Besides this,
\begin{equation}
x_n=y_{n}=z_{n}=0\quad (0<|n|< p).
                         \label{4.7}
\end{equation}
By \eqref{4.6}, $y_p=z_p=0$. Hence $x_p\neq0$ by the definition of $ p$. Using the homogeneity mentioned at the beginning of the section, we can assume without lost of generality that
\begin{equation}
|x_{ p}|=1.
                         \label{4.8}
\end{equation}

Setting $n_1=n,\ n_2= p-n$ in \eqref{1.13}, we have
$$
x_{n}y_{ p-n}-x_{ p}z_{ p-n}+y_{ p}z_{-n}=0\quad(n\neq0, p-n\neq0).
$$
The third summand on the left-hand side is equal to zero by \eqref{4.6} and the equation is simplified to the following one:
$$
x_{n}y_{ p-n}-x_{ p}z_{ p-n}=0\quad(n\neq0, p-n\neq0).
$$
If $1- p\le n\le -1$, then the first summand on the left-hand side is equal to zero by \eqref{4.7} and we obtain $x_{ p}z_{ p-n}=0\ (1- p\le n\le -1)$. Since  $x_{ p}\neq0$, we get $z_{ p-n}=0\ (1- p\le n\le -1)$ that can be written as follows: $z_{n}=0\ ( p+1\le n\le 2 p-1)$. Together with \eqref{4.6}--\eqref{4.7}, this gives
\begin{equation}
z_{n}=0\quad(0<|n|\le2 p).
                         \label{4.9}
\end{equation}

Next, we set $n_1= p,\ n_2=n$ in \eqref{1.13}
$$
x_{ p}y_{n}-x_{n+ p}z_{n}+y_{n+ p}z_{- p}=0\quad(n\neq0,n+ p\neq0).
$$
The third summand on the left-hand side is equal to zero by \eqref{4.6} and the equation is simplified to the following one:
$$
x_{ p}y_{n}-x_{n+ p}z_{n}=0\quad(n\neq0,n+ p\neq0).
$$
If $0<|n|\le2 p$, then the second summand on the left-hand side is equal to zero by \eqref{4.9} and we obtain $x_{ p}y_{n}=0\ (0<|n|\le2 p)$. Since $x_{ p}\neq0$, we conclude $y_{n}=0\ (0<|n|\le2 p)$. Together with \eqref{4.9}, this gives
\begin{equation}
y_n=z_{n}=0\quad(0<|n|\le2 p).
                         \label{4.10}
\end{equation}

Some difficulty appears in our further arguments. Let us try to take one more step.
Setting $n_1=n,\ n_2=2 p-n$ in \eqref{1.13}, we have
$$
x_{n}y_{2 p-n}-x_{2 p}z_{2 p-n}+y_{2 p}z_{-n}=0\quad(n\neq0,2 p-n\neq0).
$$
The third summand on the left-hand side is equal to zero by \eqref{4.6} and the equation is simplified to the following one:
$$
x_{n}y_{2 p-n}-x_{2 p}z_{2 p-n}=0\quad(n\neq0,2 p-n\neq0).
$$
If $1- p\le n\le -1$, then $x_n=0$ by \eqref{4.7} and we obtain
$$
x_{2 p}z_{2 p-n}=0\quad (1- p\le n\le -1).
$$
This can be written as follows:
$$
x_{2 p}z_{n}=0\quad (2 p+1\le n\le 3 p-1).
$$
By \eqref{4.6}, $z_{3p}=0$. Therefore the previous formula is strengthened to the following one:
\begin{equation}
x_{2 p}z_{n}=0\quad (2 p+1\le n\le 3 p).
                         \label{4.11}
\end{equation}

Next, we set $n_1=2 p,\ n_2=n$ in \eqref{1.13}
$$
x_{2 p}y_{n}-x_{2 p+n}z_{n}+y_{2 p+n}z_{-2 p}=0\quad(n\neq0,2 p+n\neq0).
$$
The third summand on the left-hand side is equal to zero by \eqref{4.6} and the equation is simplified to the following one:
$$
x_{2 p}y_{n}-x_{2 p+n}z_{n}=0\quad(n\neq0,2 p+n\neq0).
$$
If $1- p\le 2 p+n\le -1$, then $x_{2 p+n}=0$ by \eqref{4.7} and we obtain
$$
x_{2 p}y_{n}=0\quad(-3 p+1\le n\le -2 p-1).
$$
Changing the sign of $n$ and using $y_{-n}=\overline{y_n}$, we write this in the form
$$
x_{2 p}y_{n}=0\quad(2 p+1\le n\le 3 p-1).
$$
By \eqref{4.6}, $y_{3p}=0$. Therefore the previous formula is strengthened to the following one:
\begin{equation}
x_{2 p}y_{n}=0\quad(2 p+1\le n\le 3 p).
                         \label{4.12}
\end{equation}

Comparing \eqref{4.11} and \eqref{4.12}, we arrive to the alternative
\begin{equation}
\mbox{either}\ x_{2 p}=0\ \mbox{or}\ (y_{n}=z_n=0\ \mbox{for}\ 0<|n|\le 3 p).
                         \label{4.13}
\end{equation}
Because of the alternative, we have to consider two possible cases:

{\bf Case 1.} $x_{kp}\neq0$ for all $k>0$;

{\bf Case 2.} There exists $k\ge2$ such that $x_{kp}=0$.

We will first finish the proof in the first case. Thus, we assume that
\begin{equation}
x_{kp}\neq0\quad(k>0).
                         \label{4.14}
\end{equation}
By induction on $k$, we prove that
\begin{equation}
y_{n}=z_n=0\ \mbox{for}\ 0<|n|\le (k+1)p.
                         \label{4.15}
\end{equation}
After \eqref{4.15} had been proven for all $k$, we would have $y_{n}=z_n=0$ for all $n\neq0$, i.e., our solution is a one-dimensional one. This proves Theorem \ref{Th1.3} in the first case.

We observe that \eqref{4.15} holds for $k=1$ and for $k=2$ as is seen from \eqref{4.10} and \eqref{4.13}--\eqref{4.14}. Assume the validity of \eqref{4.15} for some $k\ge2$.

Setting $n_1=n,\ n_2=(k+1)p-n$ in \eqref{1.13}, we have
$$
x_{n}y_{(k+1)p-n}-x_{(k+1)p}z_{(k+1)p-n}+y_{(k+1)p}z_{-n}=0\quad\big(n\neq0,(k+1)p-n\neq0\big).
$$
The third summand on the left-hand side is equal to zero by \eqref{4.6} and the equation is simplified to the following one:
$$
x_{n}y_{(k+1)p-n}-x_{(k+1)p}z_{(k+1)p-n}=0\quad(n\neq0,(k+1)p-n\neq0).
$$
If $1- p\le n\le -1$, then $x_n=0$ by \eqref{4.7} and we obtain
$$
x_{(k+1)p}z_{(k+1)p-n}=0\quad (1- p\le n\le -1).
$$
This can be written as:
$$
x_{(k+1)p}z_{n}=0\quad \big((k+1)p+1\le n\le (k+2)p-1\big).
$$
Since $x_{(k+1)p}\neq0$ by the assumption \eqref{4.14}, we get
\begin{equation}
z_{n}=0\quad \big((k+1)p+1\le n\le (k+2)p-1\big).
                         \label{4.16}
\end{equation}

Next, we set $n_1=(k+1)p,\ n_2=n$ in \eqref{1.13}
$$
x_{(k+1)p}y_{n}-x_{(k+1)p+n}z_{n}+y_{(k+1)p+n}z_{-(k+1)p}=0\quad\big(n\neq0,(k+1)p+n\neq0\big).
$$
The third summand on the left-hand side is equal to zero by \eqref{4.6} and the equation is simplified to the following one:
$$
x_{(k+1)p}y_{n}-x_{(k+1)p+n}z_{n}=0\quad\big(n\neq0,(k+1)p+n\neq0\big).
$$
If $1- p\le (k+1)p+n\le -1$, then $x_{(k+1)p+n}=0$ by \eqref{4.7} and we obtain
$$
x_{(k+1)p}y_{n}=0\quad\big(-(k+2)p+1\le n\le -(k+1)p-1\big).
$$
Since $x_{(k+1)p}\neq0$ by the assumption \eqref{4.14}, we get
$$
y_{n}=0\quad\big(-(k+2)p+1\le n\le -(k+1)p-1\big).
$$
Changing the sign of $n$ and using $y_{-n}=\overline{y_n}$, we write the result in the final form
\begin{equation}
y_{n}=0\quad\big((k+1)p+1\le n\le ((k+2)p-1\big).
                         \label{4.17}
\end{equation}
Uniting \eqref{4.16} and \eqref{4.17}, we have
$$
y_{n}=z_n=0\quad\big((k+1)p+1\le n\le ((k+2)p-1\big).
$$
By \eqref{4.6}, $y_{(k+2)p}=z_{(k+2)p}=0$. Therefore the previous formula is strengthened to the following one:
$$
y_{n}=z_n=0\quad\big((k+1)p+1\le n\le ((k+2)p\big).
$$
Using the parity condition \eqref{1.14}, it can be written in the more strength form:
$$
y_{n}=z_n=0\quad\big((k+1)p+1\le |n|\le ((k+2)p\big).
$$
Together with the induction assumption \eqref{4.15}, this gives
$$
y_{n}=z_n=0\ \mbox{for}\ 0<|n|\le (k+2)p.
$$
This finishes the induction step.

We proceed to considering the second case that is characterized by the existence of $k\ge2$ such that $x_{kp}=0$. Let $k_0$ be the minimal of such $k$. Thus,
$$
x_{kp}\neq0\ (1\le k\le k_0-1);\quad x_{k_0p}=0.
$$
We can reproduce first $k_0-1$ steps of the above-presented induction proof in order to prove that
\begin{equation}
y_{n}=z_n=0\ \mbox{for}\ 0<|n|\le k_0p.
                         \label{4.19}
\end{equation}

We prove Theorem \ref{Th1.3} by contradiction. Recall that the system \eqref{1.13}--\eqref{1.14} has no two-dimensional solution. Let us assume, contrary to the assertion of Theorem \ref{Th1.3}, that the solution $(x_n,y_n,z_n)$ is three-dimensional one, i.e., each of the sequences $\{y_n\}_{n\in\Z\setminus\{0\}}$ and
$\{z_n\}_{n\in\Z\setminus\{0\}}$ has at least one nonzero term. Let $q$ be the minimal of positive integers $n$ such that $y_n\neq0$. Similarly, let $r$ be the minimal of positive integers $n$ such that $z_n\neq0$. In view of symmetries \eqref{4.1}, we can assume without lost of generality that $q\le r$. Thus,
$$
y_{n}=0\ \mbox{for}\ 0<|n|< q,\ y_q\ne0;
$$
\begin{equation}
z_{n}=0\ \mbox{for}\ 0<|n|< r,\ z_r\ne0.
                         \label{4.20}
\end{equation}
The sequence
$({\tilde x}_n,{\tilde y}_n,{\tilde z}_n)=(x_{q n},y_{q n},z_{q n})$ satisfies hypotheses of Lemma \ref{L4.1} since $\tilde y_1\neq0$. By the lemma, $({\tilde x}_n,{\tilde y}_n,{\tilde z}_n)$ is a one-dimensional solution, i.e., ${\tilde x}_n={\tilde z}_n=0\ (n\neq0)$. Thus,
$$
x_{q n}=z_{q n}=0\quad (n\in\Z\setminus\{0\}).
$$
In particular, $z_q=0$. Taking the inequality $q\le r$ into account and using \eqref{4.20}, we can now state that
\begin{equation}
z_{n}=0\ \mbox{for}\ 0<|n|\le q.
                         \label{4.21}
\end{equation}

Setting $n_1=n,\ n_2=q$ in \eqref{1.13}, we have
$$
x_{n}y_{q}-x_{n+q}z_{q}+y_{n+q}z_{-n}=0\quad(n\neq0, n+q\neq0).
$$
The second summand on the left-hand side is equal to zero by \eqref{4.21} and the equation is simplified to the following one:
$$
x_{n}y_{q}+y_{n+q}z_{-n}=0\quad(n\neq0, n+q\neq0).
$$
For $0<n\le q$, the second summand on the left-hand side is equal to zero by \eqref{4.21} and we obtain
$$
x_{n}y_{q}=0\quad (0<n\le q).
$$
Since $y_{q}\neq0$, we get
\begin{equation}
x_{n}=0\quad (0<n\le q).
                         \label{4.24}
\end{equation}
Recall that $p\le q$. Therefore \eqref{4.24} means in particular that $x_p=0$. This contradicts to the equality \eqref{4.8}. The contradiction finishes the proof of Theorem \ref{Th1.3} under the assumption that Lemma \ref{L4.1} is true.

\section{Proof of Lemma \ref{L4.1}}

By the hypothesis $|x_1|+|y_1|+|z_1|>0$, at least one of the numbers $x_1,y_1,z_1$ is not equal to zero. In view of symmetries \eqref{4.1} we can assume without lost of generality that $x_1\neq0$. Using the homogeneity mentioned at the beginning of the previous section, we can assume that
\begin{equation}
|x_1|=1.
                         \label{5.1}
\end{equation}

Let us prove the alternative
\begin{equation}
y_1=z_1=0\ \mbox{or}\ (y_1\neq0\ \mbox{and}\ z_1\neq0).
                         \label{5.2}
\end{equation}
Indeed, assume that $z_1=0$. Setting $n_1=n_2=1$ in \eqref{1.13}, we have
$$
x_1y_1-x_2z_1+y_2z_{-1}=0.
$$
The second and third summands on the left-hand side are equal to zero and we obtain $x_1y_1=0$. Since $x_1\neq0$, this implies that $y_1=0$. On the other hand, assume that $y_1=0$. Setting $n_1=2,\ n_2=-1$ in \eqref{1.13}, we have
$$
x_2y_{-1}-x_1z_{-1}+y_1z_{-2}=0.
$$
The first and third summands on the left-hand side are equal to zero and we obtain $x_1z_{-1}=0$. Since $x_1\neq0$, this implies $z_{-1}=0$. Hence $z_1=\overline{z_{-1}}=0$.

Lemma \ref{L4.1} is proved in different ways in the two cases presented by the alternative \eqref{5.2}. We first consider the case when
\begin{equation}
y_1=z_1=0.
                         \label{5.3}
\end{equation}

Set $n_1=n,n_2=1$ in \eqref{1.13}
$$
x_{n}y_{1}-x_{n+1}z_{1}+y_{n+1}z_{-n}=0\quad(n\neq0,n\neq-1).
$$
First two summands on the left-hand side are equal to zero by \eqref{5.3}. Hence $y_{n+1}z_{-n}=0$ that is equivalent to
\begin{equation}
y_{n+1}z_{n}=0\quad(n\neq0,n\neq-1).
                         \label{5.4}
\end{equation}

Next, set $n_1=n,n_2=1-n$ in \eqref{1.13}
$$
x_{n}y_{1-n}-x_{1}z_{1-n}+y_{1}z_{-n}=0\quad(n\neq0,n\neq1).
$$
The last summand on the left-hand side is equal to zero by \eqref{5.4}. Hence
$$
x_{n}y_{1-n}-x_{1}z_{1-n}=0\quad(n\neq0,n\neq1).
$$
We express $z_{1-n}$ from the last formula
$$
z_{1-n}=\frac{1}{x_{1}}\,x_{n}y_{1-n}\quad(n\neq0,n\neq1).
$$
This can be written in the form
\begin{equation}
z_{n}=\frac{1}{x_{1}}\,x_{1-n}y_{n}\quad(n\neq0,n\neq1).
                         \label{5.5}
\end{equation}

Finally, set $n_1=1,n_2=n$ in \eqref{1.13}
$$
x_{1}y_{n}-x_{n+1}z_{n}+y_{n+1}z_{-1}=0\quad(n\neq0,n\neq-1).
$$
The last summand on the left-hand side is equal to zero by \eqref{5.3}. Hence
$$
x_{1}y_{n}-x_{n+1}z_{n}=0\quad(n\neq0,n\neq-1).
$$
Express $y_{n}$ from the last formula
\begin{equation}
y_{n}=\frac{1}{x_{1}}\,x_{n+1}z_{n}\quad(n\neq0,n\neq-1).
                         \label{5.6}
\end{equation}

If $y_n=z_n=0$ for all $n\ge1$, then our solution is one-dimensional. Therefore we assume the existence of $n_0\ge2$ such that
\begin{equation}
y_{n}=z_{n}=0\quad\mbox{for}\quad0<|n|<n_0
                         \label{5.7}
\end{equation}
and $|y_{n_0}|+|z_{n_0}|>0$. Then
\begin{equation}
y_{n_0}\neq0,\quad z_{n_0}\neq0.
                         \label{5.8}
\end{equation}
Indeed, if for instance the equality $y_{n_0}=0$ was true, then setting $n=n_0$ in \eqref{5.5} we would obtain $z_{n_0}=0$. If $z_{n_0}=0$, then we use \eqref{5.6} in the same way.

Setting $n=n_0$ and then setting $n=-n_0-1$ in \eqref{5.4}, we have
$$
y_{n_0+1}z_{n_0}=0,\quad y_{-n_0}z_{-n_0-1}=0.
$$
This is equivalent to the equalities
$$
y_{n_0+1}z_{n_0}=0,\quad y_{n_0}z_{n_0+1}=0
$$
which give with the help of \eqref{5.8} 
$$
y_{n_0+1}=z_{n_0+1}=0.
$$

Let us write down the equation \eqref{4.3} for $n=n_0$
$$
x_{-m}y_{n_0}-x_{n_0-m}z_{n_0}=-y_{n_0-m}z_{m}\quad(1\le m\le n_0-1).
$$
The right-hand side is equal to zero by \eqref{5.7} and we obtain
\begin{equation}
x_{-m}y_{n_0}-x_{n_0-m}z_{n_0}=0\quad(1\le m\le n_0-1).
                         \label{5.10}
\end{equation}
In particular, for $m=1$ this gives
\begin{equation}
y_{n_0}=\frac{x_{n_0-1}}{x_{-1}}\,z_{n_0},
                         \label{5.11}
\end{equation}
and for $m=n_0-1$ \eqref{5.10} gives
$$
z_{n_0}=\frac{x_{-n_0+1}}{x_{1}}\,y_{n_0}.
$$
As follows from two last formulas, 
$$
(|x_{n_0-1}|^2-|x_1|^2)y_{n_0}=0.
$$
Since $y_{n_0}\neq0$, we obtain
\begin{equation}
|x_{n_0-1}|=|x_1|=1.
                         \label{5.13}
\end{equation}
Formulas \eqref{5.11} and \eqref{5.13} imply the important conclusion:
\begin{equation}
|y_{n_0}|=|z_{n_0}|.
                         \label{5.14}
\end{equation}

Now, we prove by induction on $k$ the following statement: for every integer $k\ge0$,
\begin{equation}
|x_{ln_0+1}|=1,\quad y_{ln_0+1}=z_{ln_0+1}=0\quad(0\le l\le k)
                         \label{5.15}
\end{equation}
and
\begin{equation}
|y_{ln_0}|=|z_{ln_0}|\neq0\quad(0< l\le k+1).
                         \label{5.16}
\end{equation}
If \eqref{5.15} was proved for all $k$, this would contradict to the decay condition $x_n\rightarrow0$.

For $k=0$, equalities \eqref{5.15}--\eqref{5.16} hold as is seen from \eqref{5.1}, \eqref{5.3}, \eqref{5.8} and \eqref{5.14}.

For the induction step, we write down equations \eqref{4.2}--\eqref{4.4} for $n=(k+1)n_0+1$
\begin{equation}
z_{(k+1)n_0-m+1}x_{(k+1)n_0+1}-z_{-m}y_{(k+1)n_0+1}=x_{m}y_{(k+1)n_0-m+1}\quad(1\le m\le (k+1)n_0),
                         \label{5.17}
\end{equation}
\begin{equation}
x_{-m}y_{(k+1)n_0+1}-x_{(k+1)n_0-m+1}z_{(k+1)n_0+1}=-y_{(k+1)n_0-m+1}z_{m}\quad(1\le m\le (k+1)n_0),
                         \label{5.18}
\end{equation}
\begin{equation}
y_{-m}x_{(k+1)n_0+1}+y_{(k+1)n_0-m+1}z_{-(k+1)n_0-1}=x_{(k+1)n_0-m+1}z_{-m}\quad(1\le m\le (k+1)n_0).
                         \label{5.19}
\end{equation}

For $m=kn_0+1$, the equation \eqref{5.19} looks as follows:
$$
y_{-kn_0-1}x_{(k+1)n_0+1}+y_{n_0}z_{-(k+1)n_0-1}=x_{n_0}z_{-(k+1)n_0-1}.
$$
By the induction hypothesis \eqref{5.15}, $y_{-kn_0-1}=z_{-kn_0-1}=0$. The first summand on the left-hand side of the last equation is equal to zero as well as the right-hand side. The equation is simplified to the following one: $y_{n_0}z_{-(k+1)n_0-1}=0$ that is equivalent to $y_{n_0}z_{(k+1)n_0+1}=0$. Observe that $y_{n_0}\neq0$ as is seen from the induction hypothesis \eqref{5.16}. Hence
\begin{equation}
z_{(k+1)n_0+1}=0.
                         \label{5.20}
\end{equation}

Let us write down the equation \eqref{4.3} for $n=2n_0$ and $m=n_0$
$$
x_{-n_0}y_{2n_0}-x_{n_0}z_{2n_0}=-y_{n_0}z_{n_0}.
$$
By the induction hypothesis \eqref{5.16}, the right-hand side of this equation is not equal to zero. Therefore
\begin{equation}
x_{n_0}\neq0.
                         \label{5.21}
\end{equation}

For $m=n_0$, the equation \eqref{5.18} looks as follows:
$$
x_{-n_0}y_{(k+1)n_0+1}-x_{kn_0+1}z_{(k+1)n_0+1}=-y_{kn_0+1}z_{n_0}.
$$
By the induction hypothesis \eqref{5.15}, $y_{kn_0+1}=0$. The right-hand side of the last equation is equal to zero. The equation is simplified to the following one:
$$
x_{-n_0}y_{(k+1)n_0+1}-x_{kn_0+1}z_{(k+1)n_0+1}=0.
$$
With the help of \eqref{5.20}, this is transformed to the form $x_{-n_0}y_{(k+1)n_0+1}=0$. Since $x_{-n_0}\neq0$ by \eqref{5.21}, we conclude
\begin{equation}
y_{(k+1)n_0+1}=0.
                         \label{5.22}
\end{equation}

For $m=1$, the equation \eqref{5.17} looks as follows:
$$
z_{(k+1)n_0}x_{(k+1)n_0+1}-z_{-1}y_{(k+1)n_0+1}=x_{1}y_{(k+1)n_0}.
$$
By the induction hypothesis \eqref{5.15}, $z_{-1}=0$. The second summand on the left-hand side of the last equation is equal to zero. The equation is simplified to the following one:
$$
z_{(k+1)n_0}x_{(k+1)n_0+1}=x_{1}y_{(k+1)n_0}.
$$
By the induction hypotheses \eqref{5.15}--\eqref{5.16}, $|x_{1}|=1,\ |y_{(k+1)n_0}|=|z_{(k+1)n_0}|\neq0$. Therefore the last equation implies
\begin{equation}
|x_{(k+1)n_0+1}|=|x_{1}|\,\frac{|y_{(k+1)n_0}|}{|z_{(k+1)n_0}|}=1.
                         \label{5.23}
\end{equation}
Together with \eqref{5.20} and \eqref{5.22}, this proves \eqref{5.15} for $k:=k+1$.

It remains to prove \eqref{5.16} for $k:=k+1$. To this end we write down the equation \eqref{4.3} for $n=(k+2)n_0$ and $m=(k+1)n_0+1$
$$
x_{-(k+1)n_0-1}y_{(k+2)n_0}-x_{n_0-1}z_{(k+2)n_0}=-y_{n_0-1}z_{(k+1)n_0+1}.
$$
By \eqref{5.7}, $y_{n_0-1}=0$. Therefore the last equation is simplified to the following one:
$$
x_{-(k+1)n_0-1}y_{(k+2)n_0}-x_{n_0-1}z_{(k+2)n_0}=0.
$$
From this
$$
|x_{(k+1)n_0+1}||y_{(k+2)n_0}|=|x_{n_0-1}||z_{(k+2)n_0}|.
$$
By \eqref{5.13}, $|x_{n_0-1}|=1$. By \eqref{5.23}, $|x_{(k+1)n_0+1}|=1$. Therefore the last formula gives
\begin{equation}
|y_{(k+2)n_0}|=|z_{(k+2)n_0}|.
                         \label{5.24}
\end{equation}

Now, we write the equation \eqref{4.3} for $n=(k+2)n_0$ and $m=n_0$
$$
x_{-n_0}y_{(k+2)n_0}-x_{(k+1)n_0}z_{(k+2)n_0}=-y_{(k+1)n_0}z_{n_0}.
$$
By the induction hypothesis \eqref{5.16}, the right-hand side of this equality is not equal to zero. Therefore the statement \eqref{5.24} can be strengthened to the following one:
$$
|y_{(k+2)n_0}|=|z_{(k+2)n_0}|\neq0.
$$
This means the validity of \eqref{5.16} for $k:=k+1$. The induction step is done. As has been mentioned before, this proves Lemma \ref{L4.1} in the case \eqref{5.3}.

\bigskip

Now, we consider the second case of the alternative \eqref{5.2} when the equality \eqref{5.1} holds and
$$
y_{1}\neq0,\quad z_{1}\neq0.
$$

Equations \eqref{4.2}--\eqref{4.4} imply the following statement.

$(\ast)$ {\it The following is true for every integer $n\geq1$. If two of three numbers $(x_n,y_n,z_n)$ are nonzero, then the third one is also nonzero.}

Indeed, let us write the equation \eqref{4.2} for $n:=2n$ and $m=n$
$$
z_{n}x_{2n}-z_{-n}y_{2n}=x_{n}y_{n},
$$
Assume that $x_ny_n\neq0$. Then $z_n\neq0$ as is seen from the last formula. Two other possible cases are considered in the same way on using \eqref{4.3}--\eqref{4.4}.

Let us prove that
\begin{equation}
x_ny_nz_n\neq0\quad\mbox{for all}\quad n>0.
                         \label{5.26}
\end{equation}
We prove by contradiction. Assume the statement to be wrong and let $n_0\ge2$ be the minimal positive integer such that $x_{n_0}y_{n_0}z_{n_0}=0$. Then
\begin{equation}
x_ny_nz_n\neq0\quad(1\le n\le n_0-1)
                         \label{5.27}
\end{equation}
and, according to the statement $(\ast)$, at least two of three numbers $(x_{n_0},y_{n_0},z_{n_0})$ are equal to zero.
We write the equation \eqref{4.2} for $n=n_0$ and $m=1$
$$
z_{n_0-1}x_{n_0}-z_{-1}y_{n_0}=x_{1}y_{n_0-1}.
$$
If $x_{n_0}=y_{n_0}=0$, then the left-hand side is equal to zero. But the right-hand side is not equal to zero by \eqref{5.27}. We have got a contradiction. Two other possible cases are considered in the same way on using \eqref{4.3}--\eqref{4.4}. Thus, \eqref{5.26} is proved.
Now, we conclude with the help of \eqref{4.5} that the differences
\begin{equation}
|x_n|^2-|y_n|^2,\quad |x_n|^2-|z_n|^2
                         \label{5.28}
\end{equation}
are independent of $n$.

By the decay condition $|x_n|+|y_n|+|z_n|\rightarrow0$, each of sequences \eqref{5.28} converges to zero as $n\rightarrow\infty$. And since the differences are independent of $n$, we conclude
\begin{equation}
|x_n|=|y_n|=|z_n|>0\quad (n=1,2,\dots).
                         \label{5.29}
\end{equation}

According to \eqref{5.29}, we represent the complex numbers $x_n,y_n,z_n$ in the trigonometric form
\begin{equation}
x_n=r_ne^{i\alpha_n},\quad y_n=r_ne^{i\beta_n},\quad z_n=r_ne^{i\gamma_n},
                         \label{5.30}
\end{equation}
where $\{r_n\}_{n=1}^\infty$ is a sequence of positive numbers converging to zero as $n\rightarrow\infty$.
Choose $n_0\ge3$ such that
\begin{equation}
r_n<1\quad\mbox{for}\quad n\ge n_0.
                         \label{5.31}
\end{equation}

Let $n\ge3$. We write two versions of the equation \eqref{4.2} corresponding to $m=1$ and $m=n-1$
$$
\begin{aligned}
z_{n-1}x_{n}-z_{-1}y_{n}&=x_{1}y_{n-1},\\
z_{1}x_{n}-z_{-n+1}y_{n}&=y_{1}x_{n-1}.
\end{aligned}
$$
It is the system of linear equations for the unknowns $(x_n,y_n)$.
The determinant of the system is equal to $|z_{1}|^2-|z_{n-1}|^2=1-r_{n-1}^2$. By \eqref{5.31}, the determinant is positive for $n>n_0$. Solving the system, we obtain
$$
x_{n}=\frac{y_{1}z_{-1}x_{n-1}-x_{1}y_{n-1}z_{-n+1}}{1-r_{n-1}^2}\quad(n>n_0).
$$
Let us increase the value of $n$ by one in order to simplify further formulas
\begin{equation}
x_{n+1}=\frac{y_{1}z_{-1}x_{n}-x_{1}y_{n}z_{-n}}{1-r_{n}^2}\quad(n\ge n_0).
                         \label{5.32}
\end{equation}

Now, using \eqref{5.30} and \eqref{5.32}, we compute
$$
\begin{aligned}
&r_{n+1}^2=|x_{n+1}|^2=\frac{1}{(1-r_{n}^2)^2}\,\big(y_{1}z_{-1}x_{n}-x_{1}y_{n}z_{-n})(y_{-1}z_{1}x_{-n}-x_{-1}y_{-n}z_{n}\big)\\
&=\frac{1}{(1-r_{n}^2)^2}\big(|y_{1}|^2|z_{1}|^2|x_{n}|^2-x_{-1}y_{1}z_{-1}x_{n}y_{-n}z_{n}
-x_{1}y_{-1}z_{1}x_{-n}y_{n}z_{-n}
+|x_{1}|^2|y_{n}|^2|z_{n}|^2\big)\\
&=\frac{1}{(1-r_{n}^2)^2}\,\Big(r_{n}^2-r_{n}^3\big(e^{i(-\alpha_1+\beta_1-\gamma_1+\alpha_{n}-\beta_{n}+\gamma_{n})}
+e^{i(\alpha_1-\beta_1+\gamma_1-\alpha_{n}+\beta_{n}-\gamma_{n})}\big)+r_{n}^4\Big).
\end{aligned}
$$
Thus,
$$
r_{n+1}=\frac{r_{n}}{1-r_{n}^2}\,\big(1-2r_{n}\cos(\alpha_{n}-\beta_{n}+\gamma_{n}-\alpha_1+\beta_1-\gamma_1)+r_{n}^2\big)^{1/2}\quad(n\ge n_0).
$$

Let us rewrite the last formula as follows:
$$
\frac{r_{k}}{r_{k-1}}=\frac{\big(1-
2r_{k-1}\cos(\alpha_{k-1}-\beta_{k-1}+\gamma_{k-1}-\alpha_1+\beta_1-\gamma_1)+r_{k-1}^2\big)^{1/2}}{1-r_{k-1}^2}\quad(k> n_0).
$$
This implies the inequality
\begin{equation}
\frac{r_{k}}{r_{k-1}}\ge\frac{1}{1+r_{k-1}}\quad(k> n_0).
                         \label{5.33}
\end{equation}

Let $n>n_0$. Taking the product of inequalities \eqref{5.33} for $k=n_0+1,n_0+2,\dots,n$, we obtain
$$
r_{n}\ge\frac{r_{n_0}}{\prod\limits_{k=n_0+1}^{n}(1+r_{k})}.
$$
Hence
\begin{equation}
r_{n}\ge\frac{r_{n_0}}{\prod\limits_{k=1}^\infty(1+r_{k})}.
                         \label{5.34}
\end{equation}

By the decay condition \eqref{1.15}, $\sum_{n=1}^\infty r_{n}<\infty$.
As well known, it is equivalent to the statement $\prod_{n=1}^\infty(1+r_{n})<\infty$.
Thus, there is a finite positive number independent of $n$ on the right-hand side of \eqref{5.34}. This contradicts to the decay condition $r_n\rightarrow0$ as $n\rightarrow\infty$. This finishes the proof of Lemma \ref{L4.1}.

\end{document}